\documentclass{amsart}

\usepackage{amsmath}
\usepackage{amssymb}
\usepackage{amsfonts}
\usepackage{latexsym}
\begin{document}


\newcommand{\pntwo}{\mathbb{P}^{n_2} \times \cdots \times \mathbb{P}^{n_k}}
\newcommand{\bj}{\underline{j}}
\newcommand{\bi}{\underline{i}}
\newcommand{\dep}{\operatorname{depth}}
\newcommand{\hit}{\operatorname{ht}}
\newcommand{\opdd}{\overline{P}_{d_1,d_2}}
\newcommand{\opdk}{\overline{P}_{d_1,\ldots,d_k}}
\newcommand{\op}{\overline{P}}
\newcommand{\xab}{X^{\underline{a}}Y^{\underline{b}}}
\newcommand{\xka}{X_1^{\underline{\alpha}_1}\cdots X_k^{\underline{\alpha}_k}}
\newcommand{\bb}{{\bf b}}
\newcommand{\ab}{(\underline{a},\underline{b})}
\newcommand{\alphak}{(\underline{\alpha}_1,\ldots,\underline{\alpha}_k)}
\newcommand{\abo}{(\overline{\underline{a},\underline{b}})}
\newcommand{\alphako}{\overline{(\underline{\alpha}_1,\ldots,\underline{\alpha}_k)}}
\newcommand{\opn}{\mathcal{O}_{\pr^n}}
\newcommand{\opthree}{\mathcal{O}_{\pr^3}}
\newcommand{\Iz}{I_{\Z}}
\newcommand{\Ixp}{I_{\xp}}
\newcommand{\Z}{\mathbb{Z}}
\newcommand{\xpi}{\X_{P_i}}
\newcommand{\xq}{\X_{Q_1}}
\newcommand{\xqi}{\X_{Q_i}}
\newcommand{\xp}{\X_{P_1}}
\newcommand{\lp}{L_{P_1}}
\newcommand{\ax}{\alpha_{\X}}
\newcommand{\bx}{\beta_{\X}}
\newcommand{\kxo}{k[x_0,\ldots,x_n]}
\newcommand{\kx}{k[x_1,\ldots,x_n]}
\newcommand{\popo}{\mathbb{P}^1 \times \mathbb{P}^1}
\newcommand{\pr}{\mathbb{P}}
\newcommand{\pn}{\mathbb{P}^n}
\newcommand{\pnpm}{\mathbb{P}^n \times \mathbb{P}^m}
\newcommand{\pnk}{\mathbb{P}^{n_1} \times \cdots \times \mathbb{P}^{n_k}}
\newcommand{\X}{\mathbb{X}}
\newcommand{\Y}{\mathbb{Y}}
\newcommand{\N}{\mathbb{N}}
\newcommand{\M}{\mathbb{M}}
\newcommand{\Q}{\mathbb{Q}}
\newcommand{\Ix}{I_{\X}}
\newcommand{\pix}{\pi_1(\X)}
\newcommand{\pixt}{\pi_2(\X)}
\newcommand{\pipi}{\pi_1^{-1}(P_i)}
\newcommand{\piqi}{\pi_2^{-1}(Q_i)}
\newcommand{\qpi}{Q_{P_i}}
\newcommand{\pqi}{P_{Q_i}}
\newcommand{\pitk}{\pi_{2,\ldots,k}}

\newtheorem{theorem}{Theorem}[section]
\newtheorem{corollary}[theorem]{Corollary}
\newtheorem{proposition}[theorem]{Proposition}
\newtheorem{lemma}[theorem]{Lemma}
\newtheorem{question}[theorem]{Question}
\newtheorem{problem}{Problem}
\newtheorem{conjecture}[theorem]{Conjecture}

\newenvironment{remark}
{\vspace{.15cm}
\refstepcounter{theorem}
\noindent{\bf Remark \thetheorem.}}
{\vspace{.15cm} }
\newenvironment{example}
{\vspace{.15cm}
\refstepcounter{theorem}
\noindent{\bf Example \thetheorem.}}
{\vspace{.15cm} }
\newenvironment{definition}
{\vspace{.15cm}
\refstepcounter{theorem}
\noindent{\bf Definition \thetheorem.}}
{\vspace{.15cm} }


\title{The Border of the Hilbert Function of a set of points
in $\pnk$}
\thanks{Updated: October 2, 2001}
\author{Adam Van Tuyl}
\address{Department of Mathematical Sciences \\ 
Lakehead University \\ 
Thunder Bay, ON P7B 5E1, Canada}
\email{avantuyl@sleet.lakeheadu.ca}

\begin{abstract}
We describe the eventual behaviour of the Hilbert function of a set
of distinct points in $\pnk$.  As a consequence of this result, we show that 
the Hilbert function of a set of points in $\pnk$ can be determined
by computing the Hilbert function at only a finite number of values.
Our result extends the result that
the Hilbert function of a set of points in $\pr^n$ stabilizes
at the cardinality of the set  of points.
Motivated by our result, we introduce the 
notion of the {\it border} of the Hilbert function of a set of points.
By using
the Gale-Ryser Theorem, a classical result about $(0,1)$-matrices,
we characterize all the possible borders for the Hilbert function
of a set of distinct points in $\popo$.
\end{abstract}

\maketitle


\section{Introduction}

The {\em Hilbert function} of a set of points in $\pr^n$ is the basis 
for many questions about sets of points.  
To any set of points, we can associate an
algebraic object which we call the {\em coordinate ring}.  The Hilbert
function is used to obtain, among other things, algebraic information about
the coordinate ring and geometric information about the set of points. 
The papers ~\cite{GeGrRo},~\cite{GHS}, ~\cite{GeMaRo}, 
~\cite{MR}, ~\cite{Mar}, and  ~\cite{L}, are  just a partial list of the 
papers that study the connection between a set of points and its Hilbert 
function. As a tool for studying sets of points, the Hilbert function is 
extremely useful due, in part, to  a result of  Geramita, Maroscia, and 
Roberts ~\cite{GeMaRo} which gives a precise description of which 
functions can be the Hilbert function of a set of points in $\pr^n$.

In this paper we wish  to extend the study of collections of points in 
$\pr^n$ to collections of points in the multi-projective space $\pnk$.  
This is an area, to our knowledge, that has seen
little exploration.  The first foray into this territory, 
of which we are aware, appears to be a series
of papers, authored by Giuffrida, Maggioni, and Ragusa
(see \cite{GuMaRa3},~\cite{GuMaRa2},~\cite{GuMaRa1}), on points that lie on 
the quadric surface $\mathcal{Q} \subseteq \pr^3$.  Because 
$\mathcal{Q} \cong \popo$, some  of the results of Giuffrida, {\it et al.} 
can be translated into results about points in multi-projective space.  
However, there remain many unanswered questions about sets of points in $\pnk$.

This paper will focus on the Hilbert functions of sets of points in 
$\pnk$.  Because the characterization of the Hilbert functions of points in 
$\pr^n$ due to Geramita, {\it et al.} \cite{GeMaRo} plays  such an important 
r\^ole in the study of those sets, a generalization of this characterization
should be a primary objective.  We state this question formally:

\begin{question} 
\label{guidingquestion}
What can be the Hilbert function of a set of points in $\pnk$?
\end{question}

If $k=1$, then, as already noted, a solution due to Geramita,
{\it et al.} exists.  If $k \geq 2$,
then the problem remains open.  One reason for the difficulty
of this question is that the associated coordinate ring is an 
example of an $\N^k$-graded ${\bf k}$-algebra, a type of ring
whose Hilbert function is not fully understood.
Some results concerning the  Hilbert function of a multi-graded 
${\bf k}$-algebra 
have been established, as is evident in 
~\cite{ACD},~\cite{Bh},~\cite{CaDeRo},~\cite{KMV},
~\cite{RV},~\cite{S},~\cite{VW1},~\cite{VW2}.
However, the question of what functions can be the 
Hilbert function of a multi-graded ring remains an open problem,
except in the case of standard graded rings.   
In this situation, i.e., rings graded in the usual sense, then 
we have {\em Macaulay's Theorem} ~\cite{Ma} which
characterizes all functions that can be the Hilbert
function of a finitely generated graded {\bf k}-algebra.
Macaulay's Theorem was used by Geramita, {\it et al.}
to classify all the possible Hilbert functions of points in $\pr^n$.

In this paper we examine a weaker version of Question ~\ref{guidingquestion}
by asking about the eventual behaviour of the Hilbert function
of a set of points in $\pnk$.  
For sets of points in $\pr^n$, the following well known result 
describes the eventual behaviour of the Hilbert function:

\begin{proposition}
\label{eventual1} 
Let $\X  \subseteq \pn$ be a collection of $s$ distinct points.  If 
$H_{\X}$ is the Hilbert function of $\X$, then 
$H_{\X}(i) = s$ for all $i \geq s-1$. 
\end{proposition}
The main result of this paper (cf. Corollary ~\ref{bordercor}) is a 
generalization of this result to  sets of points $\pnk$.
We observe that the above proposition has two consequences
for the Hilbert function of a set of points in $\pn$.  First,
to calculate $H_{\X}(i)$ for all $i \in \N$, we need
to calculate $H_{\X}(i)$ for only a finite number of $i$.  Second,
numerical information about $\X$, in this case the cardinality of 
$\X$, tells us for which $i$ we need to compute $H_{\X}(i)$
in order to determine the Hilbert function for all $i \in \N$.

The generalization of Proposition ~\ref{eventual1} for sets of 
distinct points  in $\pnk$
that we present in this paper will also have  analogous consequences.
Specifically, if $H_{\X}:\N^k \rightarrow \N$ is
the Hilbert function of $\X$, a set of distinct points, 
we demonstrate that to compute 
$H_{\X}(\underline{i})$ for all $\underline{i} \in \N^k$, we need to compute 
$H_{\X}(\underline{i})$ for only a finite number of $\underline{i} \in \N^k$.
The other values of $H_{\X}(\bi)$ are then easily determined
from our generalization of Proposition ~\ref{eventual1}.
Moreover, the $\underline{i}$ for which we need to compute
$H_{\X}(\underline{i})$ can be determined from numerical information
about the set $\X$.

Motivated by this result, we define the {\it border} of a Hilbert 
function of a set of points in $\pnk$.  
The border divides the values of the Hilbert function into two
sets:  those values which need to be computed, and those 
values which rely on our result describing the eventual behaviour
of the Hilbert function.  At the end of the paper, we specialize
to the case of sets of distinct points in $\popo$, and show
how to classify all the possible borders by using the Gale-Ryser Theorem,
a classical result about $(0,1)$-matrices.

This paper is structured as follows.  In Section 2,  we introduce  
multi-graded rings, multi-projective spaces, and Hilbert functions.
In Section 3, we give some elementary facts about the coordinate
ring associated to a set of points in $\pnk$.  Many of these results
generalize well known results about points in $\pr^n$.
In Section 4, we give the main result of this paper.
We also define  the {\it border} of a Hilbert function of a set of points in
$\pnk$.   In the final section, we restrict our focus to sets of points
in $\popo$ and their border.  This section builds upon the earlier work of
Giuffrida, {\it et al.} (\cite{GuMaRa3}, \cite{GuMaRa2},\cite{GuMaRa1}).
We begin this section by recalling some relevant facts from combinatorics 
about $(0,1)$-matrices and  partitions.

Many of the results in this paper have their genesis in examples.  
Instrumental in generating these examples was the 
commutative algebra program {\tt CoCoA} ~\cite{C}.  The results in this 
paper were part of the my Ph.D. thesis 
~\cite{VT}. 

In this paper ${\bf k}$ will denote an algebraically
closed field of characteristic zero.

\section{Multi-graded rings, multi-projective spaces, and Hilbert functions}

In this section we recall the relevant facts and definitions
about multi-graded rings, multi-projective spaces,  and their Hilbert 
functions.  Many of these results appear to be well known, although we could 
not find a standard reference for them.

Let $\N:=\{0,1,2,\ldots\}$.  If $(i_1,\ldots,i_k) \in \N^k$,
then we denote $(i_1,\ldots,i_k)$ by $\bi$.  
We set $|\bi| := \sum_h i_h$.
If $\bi,\bj \in \N^k$, then $\bi + \bj := (i_1 + j_1, \ldots,i_k + j_k)$.
We write $\bi \geq \bj$ if $i_h \geq j_h$ for every $h= 1,\ldots,k$.  
This ordering is a partial ordering on the elements of $\N^k$.  
We also observe that $\N^k$ is a semi-group generated by $\{e_1,\ldots,e_k\}$
where $e_i$ is the $i^{th}$ standard basis vector of $\N^k$,
that is, $e_i := (0,\ldots,1,\ldots,0)$ with $1$ being in the
$i^{th}$ position.

An $\N^k$-{\it graded ring} (or simply a {\it multi-graded ring} if 
$k$ is clear from the context) 
is a ring $R$ that has a direct sum
decomposition $R = {\displaystyle \bigoplus_{\bi \in \N^k} R_{\bi}}$
such that  $R_{\bi}R_{\bj} \subseteq R_{\bi + \bj}$ for all
$\bi,\bj \in \N^k$.  
We sometimes write $R_{(i_1,\ldots,i_k)}:=R_{\bi}$ as
$R_{i_1,\ldots,i_k}$ to simplify our notation.
An element $x \in R$ is said to
be $\N^k$-{\it homogeneous} (or simply {\it homogeneous}
if it is clear that $R$ is $\N^k$-graded) if $x \in R_{\bi}$
for some $\bi \in \N^k$.  If $x$ is homogeneous, then
$\deg x := \bi$.  If $k=2$, then we sometimes say that $R$
is {\it bigraded} and $x$ is {\it bihomogeneous}.

We now assume
that $R = {\bf k}[x_{1,0},\ldots,x_{1,n_1},x_{2,0},\ldots,x_{2,n_2},
\ldots,x_{k,0},\ldots,x_{k,n_k}]$.  
We induce an $\N^k$-grading on $R$
by setting $\deg x_{i,j} = e_i$.  If $k=2$, then we sometimes
write $R$ as $R = {\bf k}[x_0,\ldots,x_n,y_0,\ldots,y_m]$ with
$\deg x_i = (1,0)$ and $\deg y_i = (0,1)$.  

If $m \in R$ is a monomial,
then 
\[ m = x_{1,0}^{a_{1,0}}\cdots x_{1,n_1}^{a_{1,n_1}} 
x_{2,0}^{a_{2,0}}\cdots x_{2,n_2}^{a_{2,n_2}}\cdots
x_{k,0}^{a_{k,0}}\cdots x_{k,n_k}^{a_{k,n_k}}.\]
We denote $m$ by $X_1^{\underline{a}_1}X_2^{\underline{a}_2}\cdots
X_k^{\underline{a}_k}$ where $\underline{a}_i \in \N^{n_i+1}$.
It follows that $\deg m = (|\underline{a}_1|,|\underline{a}_2|,\ldots,
|\underline{a}_k|)$.  If $F \in R$, then we can write $F = F_1 + 
\cdots + F_r$ where each $F_i$ is homogeneous.  The $F_i$'s
are called the {\it homogeneous terms} of $F$.  

For every $\bi \in \N^k$, the set 
$R_{\bi}$ is a finite dimensional vector space over ${\bf k}$.  A basis for 
$R_{\bi}$
as a vector space is the set of monomials
\[
\left\{m =  X_1^{\underline{a}_1}X_2^{\underline{a}_2}\cdots
X_k^{\underline{a}_k} \in R \left|~\deg m =
 (|\underline{a}_1|,|\underline{a}_2|,\ldots,
|\underline{a}_k|) = \bi \right\}\right..
\]
It follows that $\dim_{\bf k} R_{\bi} = \binom{n_1+i_1}{i_1}\binom{n_2 + i_2}{i_2}
\cdots \binom{n_k + i_k}{i_k}$.

Suppose that $I = (F_1,\ldots,F_r) \subseteq R$ is an ideal.  If each
$F_j$ is $\N^k$-homogeneous, then we say $I$ is an $\N^k$-{\it homogeneous
ideal} (or simply, a {\it homogeneous ideal}).   It can be shown
that $I$ is homogeneous if and only if for every $F \in I$, all of
$F$'s homogeneous terms are in $I$.  

If $I \subseteq R$ is any ideal, then we define $I_{\bi} := I \cap R_{\bi}$
for every $\bi \in \N^k$.  It follows that each $I_{\bi}$
is a subvector space of $R_{\bi}$.  Clearly $I \supseteq {\displaystyle
\bigoplus_{\bi \in \N^k} I_{\bi}}$.  If $I$ is $\N^k$-homogeneous,
then 
$I = {\displaystyle \bigoplus_{\bi \in \N^k} I_{\bi}}$
because the homogeneous terms of $F$ belong to $I$ if $F \in I$.  

Let $I \subseteq R$ be a homogeneous ideal and consider
 the quotient ring $S = R/I$.  The ring $S$ inherits an
$\N^k$-graded ring structure if we define $S_{\bi} = (R/I)_{\bi} :=
R_{\bi}/I_{\bi}$, and hence, $S = {\displaystyle \bigoplus_{\bi \in \N^k}
(R/I)_{\bi}}$.  

Suppose that $S=R/I$ is an $\N^k$-graded ring.  The numerical function
$H_S:\N^k \rightarrow \N$ defined by
\[H_S(\bi) := \dim_{\bf k} (R/I)_{\bi} = \dim_{\bf k} R_{\bi} - \dim_{\bf k} I_{\bi}\]
is the {\it Hilbert function of $S$}. 

\begin{remark}
If $k = 1$, then a precise description of what functions
can be the Hilbert function of a standard graded ${\bf k}$-algebra $S$,
i.e., $S = R/I$ for some ideal $I \subseteq R$,
was first given by Macaulay ~\cite{Ma}.  If $k \geq 2$,
then it remains an open problem to give such a description.  Some necessary
conditions for the Hilbert function of an $\N^2$-graded ${\bf k}$-algebra
were given by Aramova, {\it et al.} ~\cite{ACD}.
\end{remark}

We now extend the classical definition of projective space to 
multi-projective space.
We define the {\it multi-projective space $\pnk$} to be
\[\pnk := \left\{
\begin{tabular}{l}
$((\underline{a}_1),\ldots,(\underline{a}_k))\in
{\bf k}^{n_1+1}\times \cdots \times {\bf k}^{n_k +1}$ \\
\text{with no $\underline{a}_i = (a_{i,0},\ldots,a_{i,n_i}) = \underline{0}$}
\end{tabular}\right\}_{\displaystyle /\sim} \]
where $(\underline{a}_1,\ldots,\underline{a}_k) \sim 
(\underline{b}_1,\ldots,\underline{b}_k)$ if there exists non-zero
$\lambda_1,\ldots,\lambda_k \in {\bf k}$ such that for all $i =1,\ldots,k$
\[\underline{b}_i = (b_{i,0},\ldots,b_{i,n_i}) = (\lambda_i a_{i,0},\ldots,
\lambda_i a_{i,n_i}) ~~\text{where}~~ \underline{a}_i = 
(a_{i,0},\ldots,a_{i,n_i}).\] 
An element of $\pnk$ is called a {\it point}.  We sometimes denote
the equivalence class of $((a_{1,0},\ldots,a_{1,n_1}),\ldots,
(a_{k,0},\ldots,a_{k,n_k}))$ by $[a_{1,0}:\cdots:a_{1,n_1}]
\times\cdots \times
[a_{k,0}:\cdots:a_{k,n_k}]$.  It follows that $[a_{i,0}:\cdots:a_{i,n_i}]$
is a point of $\pr^{n_i}$ for every $i$.

If $F \in R = {\bf k}[x_{1,0},\ldots,x_{1,n_1},\ldots,x_{k,0},\ldots,x_{n_k}]$ 
is an $\N^k$-homogeneous element of 
degree $(d_1,\ldots,d_k)$
and $P = [a_{1,0}:\cdots:a_{1,n_1}] \times \cdots \times
[a_{k,0}:\cdots:a_{k,n_k}]$ is a point of $\pnk$, then 
\[F(\lambda_1a_{1,0},\ldots,\lambda_2 a_{2,0},\ldots,\lambda_k a_{k,0},\ldots)
= \lambda_1^{d_1}\lambda_2^{d_2}\cdots\lambda_k^{d_k} 
F(a_{1,0},\ldots,a_{2,0},\ldots,a_{k,0},\ldots).\]
To say that $F$ vanishes at a point of $\pnk$ is, therefore, a 
well-defined notion.

If $T$ is any collection of $\N^k$-homogeneous elements of $R$, then
we define
\[{\bf V}(T) := \{P \in \pnk ~~|~~ F(P) = 0 ~\text{for all}~ F \in T\}.\]
If $I$ is an $\N^k$-homogeneous ideal of $R$, then ${\bf V}(I) =
{\bf V}(T)$ where $T$ is the set of all homogeneous elements of $I$.  
If $I = (F_1,\ldots,F_r)$, then 
${\bf V}(I) = {\bf V}(F_1,\ldots,F_r)$.

The multi-projective space $\pnk$ can be endowed with a topology
by defining the {\it closed sets} to be all subsets of $\pnk$ of the
form ${\bf V}(T)$ where $T$ is a collection of $\N^k$-homogeneous
elements of $R$.  If $Y$ is a subset of $\pnk$ that is closed
and irreducible with respect to this topology, then we say $Y$
is a {\it multi-projective variety}, or simply, a {\it variety}.

If $Y$ is any subset of $\pnk$, then we set 
\[{\bf I}(Y) := \{F \in R~|~ F(P) = 0 ~\text{for all}~ P \in Y\}.\]
The set ${\bf I}(Y)$ is an $\N^k$-homogeneous ideal of $R$.
We call ${\bf I}(Y)$ the {\it $\N^k$-homogeneous ideal associated to $Y$}, or
simply, the {\it ideal associated to $Y$}.
If $Y \subseteq \pnk$, then we set 
$I_Y := {\bf I}(Y)$, and we call $R/I_Y$ the 
{\it $\N^k$-homogeneous coordinate ring of $Y$}, or simply, the
{\it coordinate ring of $Y$}.  If $H_{R/I_Y}$ is the Hilbert function
of $R/I_Y$, then we sometimes write $H_{Y}$ for
$H_{R/I_Y}$, and we say $H_Y$ is the {\it Hilbert function of $Y$}. 

By adopting the proofs of the homogeneous case, we have

\begin{proposition}
\label{basicprop}

\begin{enumerate}
\item[$(i)$]  If $I_1 \subseteq I_2$ are $\N^k$-homogeneous ideals,
then ${\bf V}(I_1) \supseteq {\bf V}(I_2)$.
\item[$(ii)$] If $Y_1 \subseteq Y_2$ are subsets of $\pnk$, then
${\bf I}(Y_1) \supseteq {\bf I}(Y_2)$.
 \item[$(iii)$]  For any two subsets $Y_1,Y_2$ of 
$\pnk$, ${\bf I}(Y_1 \cup Y_2) = {\bf I}(Y_1) \cap {\bf I}(Y_2)$.
\end{enumerate}
\end{proposition}

There is also an $\N^k$-graded analog of the Nullstellensatz.
Again, the proof follows as in the graded case.

\begin{theorem}
\label{nulls}
{\em ($\N^k$-homogeneous Nullstellensatz)}  If $I \subseteq R$
is an $\N^k$-homogeneous ideal and $F \in R$ is an
$\N^k$-homogeneous polynomial with $\deg F > \underline{0}$
such that $F(P) = 0$ for all $P \in {\bf V}(I) \subseteq \pnk$,
then $F^t \in I$ some $t > 0$.
\end{theorem}

Set ${\bf m}_i := (x_{i,0},x_{i,1},\ldots,x_{i,n_i})$ for $i = 1,\ldots,k$.
An $\N^k$-homogeneous ideal $I$ of $R$ is called {\it projectively irrelevant}
if ${\bf m}_i^a \subseteq I$ for some $i \in\{1,\ldots,k\}$ and some
positive integer $a$.  An ideal $I \subseteq R$ is {\it projectively 
relevant} if it is not projectively irrelevant.  By employing
the $\N^k$-homogeneous Nullstellensatz, it can be shown that
there is a one-to-one correspondence
between the non-empty closed subsets of $\pnk$ and the $\N^k$-homogeneous
ideals of $R$ that are radical and projectively relevant.  The 
correspondence is given by $Y \mapsto {\bf I}(Y)$ and $I \mapsto
{\bf V}(I)$.  This is analogous to the well known graded case.  
For the case $k = 2$, this correspondence can be found
in Van der Waerden ~\cite{VW1},\cite{VW2}.  Van der Waerden also
asserts that for arbitrary $k$ the results are analogous to the case $k=2$.

\begin{remark}
Our construction of $\pnk$ and its subsets follows the classical
definition of the projective space $\pr^n$ as described, for example, in
~\cite[Chap. 1]{H1}.   Van der Waerden ~\cite{VW1} gives a construction
similar to the approach we have given above.
The multi-projective
space $\pnk$ can also be constructed via the modern methods of schemes.
However, since we are interested in studying sets of distinct points,
which are examples of a reduced schemes, 
the classical approach is equivalent to the schematic approach.  Hence,
we will not invoke the language of schemes.
\end{remark}


\section{The coordinate ring associated to a set of points in $\pnk$}

In this section we investigate the structure of the coordinate
ring of a set of points in $\pnk$.  We also give
some elementary properties about the Hilbert function of the
coordinate ring.   We will only consider sets of distinct points.  
Many of these results generalize well known results about sets of
points in $\pr^n$.

Let $R = {\bf k}[x_{1,0},\ldots,x_{1,n_1},\ldots,x_{k,0},\ldots,x_{k,n_k}]$
and induce an $\N^k$-grading on $R$ by setting $\deg x_{i,j} = e_i$.
\begin{proposition}
\label{pointsprop1}
For any point $P \in \pnk$, let $I_P$ be the ideal of $R$ associated
to the point $P$.  Then
\begin{enumerate}
\item[$(i)$]  $I_{P}$ is a prime ideal.
\item[$(ii)$] $I_{P} = (L_{1,1},\ldots,L_{1,n_1},L_{2,1},\ldots,L_{2,n_2},
\ldots,L_{k,1},\ldots,L_{k,n_k})$ where $\deg L_{i,j} = e_i$.
\end{enumerate}
\end{proposition}
\begin{proof}
$(i)$  If $FG \in I_{P}$, then $(FG)(P) = F(P)G(P) = 0$.  Hence,
either $F(P) = 0$ or $G(P)=0$, i.e., either $F \in I_P$ or $G \in I_P$.

$(ii)$  Suppose that $P = [a_{1,0}:\cdots:a_{1,n_1}] 
\times \cdots \times [a_{k,0}:\cdots:a_{k,n_k}]$.
For each $i \in \{1,\ldots,k\}$ there exists  $a_{i,j} \neq 0$.
Assume for the moment that $a_{i,n_i} \neq 0$ for all $i$.  We
can then assume that
$P = [a_{1,0}:\cdots:a_{1,n_1-1}:1] \times
[a_{2,0}:\cdots:a_{2,n_2-1}:1] \times \cdots \times 
[a_{k,0}:\cdots:a_{k,n_k-1}:1].$
Set 
\[I := 
\left(\begin{tabular}{l}
$x_{1,0}- a_{1,0}x_{1,n_1}, ~x_{1,1}- a_{1,1}x_{1,n_1},\ldots, 
x_{1,n_1-1}- a_{1,n_1-1}x_{1,n_1},$ \\ 
$x_{2,0}- a_{2,0}x_{2,n_2}, x_{2,1}- a_{2,1}x_{2,n_2},\ldots, 
x_{2,n_2-1}- a_{2,n_2-1}x_{2,n_2}$, \\
\hspace{.8cm}\vdots \\
$x_{k,0}- a_{k,0}x_{k,n_k}, x_{k,1}- a_{k,1}x_{k,n_k},\ldots, 
x_{k,n_k-1}- a_{k,n_k-1}x_{k,n_k}$ \\
\end{tabular}\right).\]
Then $I \subseteq I_P$ because all of the generators of $I$ vanish
at $P$.  

If we show that $I_P \subseteq I$, then we will be finished
because $\deg (x_{i,j}- a_{i,j}x_{i,n_i}) = e_i$.  To do this, we first note
that the generators of $I$ are, in fact, a Groebner basis for $I$.
By using this fact, we can show that $I$ is a prime ideal.  Indeed,
suppose that $F,G \not\in I$.  Since
$F,G \not\in I$, the division of $F$ and $G$ by the generators of $I$ 
yields
$F = F' + F''$ and $G= G' + G''$
where $F',G' \in I$ and $F'',G'' \not\in I$.  
Since the generators of $I$ are a Groebner basis,
$F'',G''$ must be polynomials in the indeterminates 
$x_{1,n_1},x_{2,n_2},\ldots,x_{k,n_k}$ alone.  
If $FG = F'G' + F''G' + F'G'' + F''G'' \in I$, then this would imply 
that $F''G'' \in I$.  But the leading term of $F''G''$ is
not in the leading term ideal of $I$, contradicting the fact that the
generators of $I$ are a Groebner basis.  
So $FG \not\in I$ and hence, $I$ is prime.

We now demonstrate that $I_P \subseteq I$.  Let $F \in I_P$.  Because
${\bf V}(I) = {\bf V}(I_{P}) = P$, the Nullstellensatz
(Theorem ~\ref{nulls}) implies that 
$F^t \in I$ for some positive integer $t$. 
But since $I$ is prime,  $F \in I$ as desired.

To complete the proof of $(ii)$, if $a_{i,n_i} =0$, then there
exists an integer $0 \leq j <n_i$ such that $a_{i,j}\neq 0$.  We then repeat
the above argument, but use $x_{i,j}$ instead of $x_{i,n_i}$ to
form the generators of $I$, and use a monomial ordering so 
that $x_{r,s} > x_{i,j}$ if $r > i$ and if $r =i$, then
$x_{i,s} > x_{i,j}$ for all $s \in \{0,\ldots,\hat{j},\ldots,n_i\}$.
\end{proof}

For each $i \in \{1,\ldots,k\}$, we define the projective morphism
$\pi_i:\pnk \rightarrow \pr^{n_i}$ by
\[[a_{1,0}:\cdots:a_{1,n_1}]\times \cdots \times
[a_{i,0}:\cdots:a_{i,n_i}]\times \cdots \times
[a_{k,0}:\cdots:a_{k,n_k}] \longmapsto [a_{i,0}:\cdots:a_{i,n_i}].\]
If $\X$ is a finite collection of distinct points in $\pnk$,
then $\pi_i(\X) \subseteq \pr^{n_i}$ is the finite set of
distinct $i^{th}$ coordinates that appear in $\X$.  The
Hilbert function of $\pi_i(\X)$ can be read from the Hilbert function
of $\X$ as we show below.

\begin{proposition}
\label{sidehilbertfcn}
Suppose that $\X \subseteq \pnk$ is a finite set of 
points with Hilbert function $H_{\X} := H_{R/\Ix}$.  Fix an
integer $i \in \{1,\ldots,k\}$.  Then the sequence
$H = \{h_j\}_{j\in \N}$, where $h_j := H_{\X}(0,\ldots,j,\ldots,0)$
with $j$ in the $i^{th}$ position, is the Hilbert function
of $\pi_i(\X) \subseteq \pr^{n_i}$.
\end{proposition}

\begin{proof}
We will prove the statement for the case $i =1$.  The
other cases follow similarly.  Let
$I = {\bf I}(\pi_1(\X)) \subseteq S ={\bf k}[x_{1,0},\ldots,x_{1,n_1}]$.
We wish to show that $(R/\Ix)_{j,0\ldots,0} \cong 
(S/I)_j$ for all $j \in \N$.
Since $R_{j,0,\ldots,0} \cong S_j$ for all $j \in \N$, it is enough
to show that $(\Ix)_{j,0,\ldots,0} \cong I_j$ for all $j \in \N$.

If $P$ is a point of $\X \subseteq \pnk$, then,  
by Proposition ~\ref{pointsprop1}, the ideal associated to $P$
is  $I_P= (L_{1,1},\ldots,L_{1,n_1},L_{2,1},\ldots,L_{2,n_2},\ldots,L_{k,1},
\ldots,L_{k,n_k})$ where $\deg L_{i,j} = e_i$.  
Let $P'$ denote $\pi_1(P) \in \pr^{n_1}$.  Then the ideal associated
to $P'$ in $S$ is $I_{P'} = (L_{1,1},\ldots,L_{1,n_1})$
where we consider $L_{1,1},\ldots,L_{1,n_1}$ as 
$\N^1$-graded elements of $S$.  There is then
an isomorphism of vector spaces
$(I_P)_{j,0,\ldots,0} = (L_{1,1},\ldots,L_{1,n_1})_{j,0,\ldots,0} \cong (I_{P'})_j$ for each positive integer $j$.

Thus, if $\X = \{P_1,\ldots,P_s\}$, then $\pi_1(\X) = \{\pi_1(P_1),\ldots,
\pi_1(P_s)\}$, and hence
\[(\Ix)_{j,0,\ldots,0} = \bigcap_{i=1}^s  (I_{P_i})_{j,0,\ldots,0} 
\cong \bigcap_{l=1}^s \left(I_{\pi_1(P_i)}\right)_j = I_j
\hspace{.5cm}\text{for all $j \in \N$.}\]
\end{proof}

The previous theorem places a necessary condition on the Hilbert
function of a set of points $\X \subseteq \pnk$.  Specifically,
certain sub-sequences of the function $H_{\X}$ must grow like the Hilbert
function of a set of points in $\pr^n$. We end this section by giving 
some more necessary conditions on the Hilbert function of a set of 
points $\X$ in $\pnk$. We will first require the following lemma.

\begin{lemma}
\label{firstnzd1}
Suppose $\X \subseteq \pnk$ is a finite set  of distinct points.  
Fix an $i \in \{1,\ldots,k\}.$  Then there exists a form $L \in R$ of 
degree $e_i$ such that $\overline{L}$  is a non-zero divisor in $R/\Ix$.
\end{lemma}

\begin{proof}
We will show only the case $i =1$.
The primary decomposition of $\Ix$ is $\Ix = \wp_1 \cap \cdots \cap \wp_s$
where $\wp_i$ is an $\N^k$-homogeneous prime ideal associated
to a point of $\X$.  The set of zero divisors of $R/\Ix$,
denoted ${\bf Z}(R/\Ix)$, are precisely the elements of
${\bf Z}(R/\Ix) = \bigcup_{i=1}^s \overline{\wp}_i.$
We want to show ${\bf Z}(R/\Ix)_{e_1} \subsetneq (R/\Ix)_{e_1}$,
or equivalently, ${\displaystyle \bigcup_{i=1}^s (\wp_i)_{e_1} 
\subsetneq R_{e_1}}$.  By Proposition ~\ref{pointsprop1} it is clear that 
$(\wp_i)_{e_1} \subsetneq R_{e_1}$ for each $i = 1,\ldots,s$.
Because the field ${\bf k}$ is infinite, the vector space $R_{e_1}$ 
cannot be expressed as a finite union of vector spaces, and hence,  
$\bigcup_{i=1}^s (\wp_i)_{e_i} \subsetneq R_{e_1}$.
\end{proof}

\begin{remark}
The above lemma implies that $\operatorname{depth} R/\Ix \geq 1$
for all sets of points $\X \subseteq \pnk$.  The depth of $R/\Ix$
will be explored more thoroughly in a future paper.
\end{remark}

\begin{proposition}
\label{simpleboundsonhx}
Let $\X$ be a set of distinct points in $\pnk$ and suppose that
$H_{\X}$ is the Hilbert function of $\X$.  
\begin{enumerate}
\item[$(i)$] For all $\bi = (i_1,\ldots,i_k) \in \N^k$, 
$H_{\X}(\bi) \leq H_{\X}(\bi + e_j)$ for all $j = 1,\ldots,k$ 
\item[$(ii)$] If $H_{\X}(\bi)
=H_{\X}(\bi + e_j),$  for some $j \in \{1,\ldots,k\}$, then 
$H_{\X}(\bi + e_j)
=H_{\X}(\bi + 2e_j).$
\end{enumerate}
\end{proposition}
\begin{proof}
$(i)$ We will only demonstrate that $H_{\X}(\bi) \leq H_{\X}(\bi + e_1) =
H_{\X}(i_1+1,i_2,\ldots,i_k)$ since the other cases follow similarly.  
By Lemma ~\ref{firstnzd1} there exists a form $L \in R$ such
that $\deg L = e_1$ and $\overline{L}$ is a non-zero divisor
is $R/\Ix$.  Hence, for
any $\bi \in \N^k$, the multiplication map
$\left( R/\Ix\right)_{\bi} \stackrel{\times \overline{L}}
{\longrightarrow} \left( R/\Ix\right)_{\bi+e_1}$
is an injective map of vector spaces.  Therefore
$H_{\X}(\bi)
= \dim_{\bf k} (R/\Ix)_{\bi} \leq  \dim_{\bf k}(R/\Ix)_{(i_1+1,i_2,\ldots,i_k)}
= H_{\X}(\bi + e_1).$

$(ii)$ We will only consider the case that $j=1$ since the other cases
are proved similarly.  By Lemma ~\ref{firstnzd1}
there exists a form $L \in R$ such that  $\deg L = e_1$
and $\overline{L}$ is a non-zero divisor in $R/\Ix$.  Thus,
for each $\bi = (i_1,\ldots,i_k) \in \N^k$, we have 
the following short exact sequence of vector spaces:
\[0 \longrightarrow  \left( R/\Ix\right)_{\bi} 
\stackrel{\times \overline{L}}{\longrightarrow} 
\left( R/\Ix\right)_{\bi+e_1} \longrightarrow \left(R/(\Ix,L)\right)_{\bi+e_1}
\longrightarrow 0.\]
If $H_{\X}(\bi) = H_{\X}(\bi +e_1)$, then this implies that
the morphism $\times \overline{L}$ is an isomorphism of vector
spaces, and thus, $(R/(\Ix,L))_{\bi +e_1} = 0$.   So
$(R/(\Ix,L))_{\bi +2e_1} = 0$ as well.  Hence, from the
short exact sequence
\[0 \longrightarrow  \left( R/\Ix\right)_{\bi+e_1} 
\stackrel{\times \overline{L}}{\longrightarrow} 
\left( R/\Ix\right)_{\bi+2e_1} \longrightarrow 
\left(R/(\Ix,L)\right)_{\bi+2e_1} \longrightarrow 0\]
we deduce that $(R/\Ix)_{\bi+e_1} \cong (R/\Ix)_{\bi + 2e_1}$.
\end{proof}

\begin{remark}
Statement $(ii)$ of the above proposition is a generalization of 
a result for points in $\pr^n$ found in  
Geramita and Maroscia (cf. Proposition 1.1 (2) of ~\cite{GeMa}).
\end{remark}


\section{The Border of the Hilbert Function for Points in $\pnk$}

In this section we present our main result which
generalizes the following well known result (see, for example,  the 
discussion before Proposition 1.3 in \cite{GeMa}) for sets of points in 
$\pn$ to sets of points in $\pnk$.
 
\begin{proposition} 
\label{eventual}
Let $\X  \subseteq \pn$ be a collection of 
$s$ distinct points.  If 
$H_{\X}$ is the Hilbert function of $\X$, then 
$H_{\X}(i) = s$ for all $i \geq s-1$. 
\end{proposition}

So, suppose $\X \subseteq \pnk$ is a collection of $s$ distinct points.  Let
$\Ix$ denote the $\N^k$-homogeneous ideal associated to $\X$ in
the $\N^k$-graded ring 
$R = {\bf k}[x_{1,0},\ldots,x_{1,n_1},\ldots,x_{k,0},\ldots,x_{k,n_k}]$
where $\deg x_{i,j} = e_i$, the $i^{th}$ standard basis vector of $\N^k$.

Let $\pi_1 : \pnk \rightarrow
\pr^{n_1}$ be the projection morphism.  The image of $\pix$ in $\pr^{n_1}$ 
is a collection of $t_1:=|\pi_1(\X)| \leq s$ points.  The set of points
$\pix$ is the set of distinct first coordinates that appear in $\X$.
For every $P_i \in \pix$, we have
\[\pipi = \left\{P_i \times Q_{i,1}, \ldots, P_i \times Q_{i,{\alpha_i}}
\right\} 
\subseteq \X\]
where $Q_{i,j} \in \pntwo$.
Set $\alpha_i:=|\pipi| \geq 1$ for all $P_i \in \pix$.  
Note that the sets $\pipi$ partition $\X$.
Let $\pitk:\pnk \rightarrow \pntwo$ be the projection morphism.  For each
$P_i \in \pix$, the set
\[\qpi := \pi_{2,\ldots,k}(\pipi) = 
\left\{Q_{i,j} ~\left|~ P_i \times Q_{i,j} \in \pipi \right\}\right.\]
is a collection of $\alpha_i$ distinct points in $\pntwo$.

If $\bj = (j_1,j_2,\ldots,j_k) \in \N^k$, then we sometimes
write $\bj$ as $(j_1,\bj^{'})$ where 
$\bj^{'} = (j_2,\ldots,j_k) \in \N^{k-1}$.
Also, recall that we write $R_{j_1,\ldots,j_k}$ for
$R_{(j_1,\ldots,j_k)}$.  If $\bj = (j_1,\bj^{'})$,
then we denote $R_{(j_1,\bj^{'})} = R_{\bj}$ by
$R_{j_1,\bj^{'}}$.
With the above notation, we have 

\begin{proposition}
\label{mainresult3}
Let $\X$ be
a set of $s$ distinct points  in $\pnk$ with $k \geq 2$,
and suppose that $\pix = \{P_1,\ldots,P_{t_1}\}$
is the set of $t_1 \leq s$ distinct first coordinates in $\X$.  
Fix any tuple $\bj =(j_2,\ldots,j_k)\in \N^{k-1}$.  
Then, for all integers
$l \geq t_1-1 = |\pix| - 1,$
\begin{eqnarray*}
\dim_{\bf k} (R/I_{\X})_{l,\bj} & = & 
\sum_{P_i \in \pix} H_{\qpi}(\bj)
\end{eqnarray*}
where $H_{\qpi}$ is the Hilbert function of the set of points
$\qpi \subseteq \pntwo$.
\end{proposition}

To prove this proposition we require the following two 
results.

\begin{proposition}
\label{rank=hilbertforpnk}
Let $\X = \{P_1,\ldots,P_s\} \subseteq \pnk$ be a set
of $s$ distinct points.
For any $\bj = (j_1,\ldots,j_k) \in \N^{k}$, let
$\{m_1,\ldots,m_N\}$ be
the $N = \binom{n_1+j_1}{j_1}\binom{n_2+j_2}{j_2}\cdots
\binom{n_k + j_k}{j_k}$ monomials of $R$ of degree $\bj$.  Set
\[M_{\bj} = \bmatrix
m_1(P_1) & \cdots & m_N(P_1) \\
\vdots &  &\vdots \\
m_1(P_s) & \cdots & m_N(P_s) \\
\endbmatrix.\]
Then $\operatorname{rk}M_{\bj} = H_{\X}(\bj)$ where $H_{\X}$ is the 
Hilbert function
of $R/I_{\X}$.
\end{proposition}

\begin{proof}
To compute $H_{\X}(\bj)$, we need to determine the number of linearly 
independent forms of degree $\bj$ that pass through $\X$.  An element
of $R$ of degree $\bj$ has the form
$F = c_1m_1 + \cdots + c_Nm_N$
where $c_i \in {\bf k}$.
If $F(P_i) = 0$, we get a linear relation among the $c_i$'s, namely,
$c_1m_1(P_i) + \cdots + c_Nm_N(P_i) = 0.$
The elements of $(I_{\X})_{\bj}$ are given by the solutions
of the system of linear equations
$F(P_i) = \cdots = F(P_s) = 0$.  We can 
rewrite this system
of equations as 
\[\bmatrix
m_1(P_1) & \cdots & m_N(P_1) \\
\vdots &  &\vdots \\
m_1(P_s) & \cdots & m_N(P_s) 
\endbmatrix
\bmatrix
c_1 \\
\vdots \\
c_N 
\endbmatrix
=
\bmatrix
0 \\
\vdots \\
0
\endbmatrix.
\]
The matrix on the left 
is $M_{\bj}$.  Now the number of linearly independent solutions is
equal to $\dim_{\bf k}
(I_{\X})_{\bj}$, and hence,
\[
\dim_{\bf k} (I_{\X})_{\bj} 
 = \#\mbox{columns of $M_{\bj}$} - \operatorname{rk} M_{\bj}  =  
N - \operatorname{rk} M_{\bj}. 
\]
Since $\dim_{\bf k} R_{\bj} = N$, we have $H_{\X}({\bj}) = 
\operatorname{rk} M_{\bj}$, as desired.
\end{proof}

\begin{proposition}
\label{seperators3}
Let $\X = \{P_1,\ldots,P_s\} \subseteq \pnk$
and suppose that
$H_{\X}(\bj) = h$.  Then we can find a subset $\X' \subseteq \X$
of $h$ elements, say $\X' = \{P_1,\ldots,P_h\}$ (after a possible reordering),
such that there exists $h$ forms $G_1,\ldots,G_h$ of degree $\bj$
with the property that 
for every $1 \leq l \leq h$,  $G_i(P_l) = 0$ if $i \neq l$, 
and $G_i(P_i) \neq 0$.
\end{proposition}

\begin{proof}
Let $\{m_1,\ldots,m_N\}$ be the 
$N = \binom{n_1+j_1}{j_1}\cdots \binom{n_k + j_k}{j_k}$ 
monomials of degree $\bj$ in $R$.  By Proposition
~\ref{rank=hilbertforpnk} the matrix
\[M_{\bj} = \bmatrix
m_1(P_1) & \cdots & m_N(P_1) \\
\vdots &  &\vdots \\
m_1(P_s) & \cdots & m_N(P_s) \\
\endbmatrix\]
has $\operatorname{rk} M_{\bj} = H_{\X}(\bj) = h$.  Without
loss of generality, we can assume that the first $h$ rows are linearly
independent.  So, let $\X' =\{P_1,\ldots,P_h\} \subseteq \X$, and let
\[
M'_{\bj} = \bmatrix
m_1(P_1) & \cdots & m_N(P_1) \\
\vdots &  &\vdots \\
m_1(P_h) & \cdots & m_N(P_h) \\
\endbmatrix.
\]
Fix an $i \in \{1,\ldots,h\}$ and let
$\X'_i =\{P_1,\ldots,\widehat{P}_i,\ldots,P_h\}$.
If we remove the $i^{th}$ row of $M'_{\bj}$, then the
rank of the resulting
matrix decreases by one.  Since the rank of the new
matrix is equal to the Hilbert function of $\X'_i$, it
follows that $\dim_{\bf k} (I_{\X'})_{\bj} + 1 = \dim_{\bf k} 
(I_{\X'_i})_{\bj}$.
Thus, there exists an element $G_i \in (I_{\X'_i})_{\bj}$ such
that $G_i$ passes through the points of $\X'_i$ but not through $P_i$.
We repeat this argument for each $i\in\{1,\ldots,h\}$
to get the desired forms.
\end{proof}

\begin{corollary}
\label{seperators}
Let $\X = \{P_1,\ldots,P_s\} \subseteq \pr^n$
be a set of $s$ distinct points. 
Then there exists $s$ forms $F_1,\ldots,F_s$ of degree $s-1$
with the property that 
for every $1 \leq l \leq s$,  $F_i(P_l) = 0$ if $i \neq l$, 
and $F_i(P_i) \neq 0$.
\end{corollary}

\begin{proof}
By Proposition ~\ref{eventual}, $H_{\X}(s-1) = s$.  Now apply
the above theorem.
\end{proof}

\begin{proof} (of Proposition ~\ref{mainresult3})
Fix a $\bj =(j_2,\ldots,j_k) \in \N^{k-1}$, let $N = N(\bj):=
\binom{n_2+j_2}{j_2}\cdots\binom{n_k + j_k}{j_k}$, and set
\begin{eqnarray*} 
 (*) & = &  \sum_{P_i \in \pix} H_{\qpi}(\bj).
\end{eqnarray*}
We will first show that $\dim_{\bf k} (R/I_{\X})_{l,\bj} \leq (*)$
for all $l \in \N$.  
Let
$\{X_1,\ldots,X_{\binom{n_1+l}{l}}\}$ be all the monomials
of degree $(l,\underline{0})$ in $R$ and let $\{Y_1,\ldots,Y_N\}$
be the $N$ monomials of degree $(0,\bj)$ in $R$.  For any
$l \in \N$, a general form $L \in R_{l,\bj}$ looks like
\begin{eqnarray*}
L & = & \left(c_{1,1}X_1 + \cdots + 
c_{1,\binom{n_1+l}{l}}X_{\binom{n_1+l}{l}}\right)Y_1  + 
\left(c_{2,1}X_1 + \cdots + c_{2,\binom{n_1+l}{l}}X_{\binom{n_1+l}{l}}
\right)Y_2 \\
  &   & + \cdots + 
\left(c_{N,1}X_1 + \cdots + c_{N,
\binom{n_1+l}{l}}X_{\binom{n_1+l}{l}}\right)Y_{N}
\end{eqnarray*}
with coefficients $c_{i,j} \in {\bf k}$.  By setting
$A_i := c_{i,1}X_1 + \cdots + c_{i,\binom{n_1+l}{l}}X_{\binom{n_1+l}{l}}$
for $i = 1,\ldots,N$, we can rewrite $L$ as 
$L = A_1Y_1 + A_2Y_2 + \cdots + A_{N}Y_{N}.$

\noindent
{\it Claim. }  Each subset $\pipi \subseteq \X$ puts at most 
$H_{\qpi}(\bj)$ linear restrictions on the forms of $R_{l,\bj}$ that 
pass through $\X$.

\noindent
{\it Proof of the Claim.} Suppose $\pipi = \{P_i \times Q_{i,1},\ldots, P_i 
\times Q_{i,{\alpha_i}}\} \subseteq \X$, and hence, the set 
 $\qpi = \{Q_{i,1},\ldots,
Q_{i,{\alpha_i}}\} \subseteq \pntwo$.  
If $L \in R_{l,j}$ vanishes at the $s$ points of $\X$,
then it vanishes on $\pipi$, and thus 
\begin{eqnarray*}
L(P_i \times Q_{i,1}) & = & A_1(P_i)Y_1(Q_{i,1}) + \cdots + 
A_{N}(P_i)Y_{N}(Q_{i,1}) = 0 \\
& \vdots & \\
L(P_i \times Q_{i,{\alpha_i}}) & = & A_1(P_i)Y_1(Q_{i,{\alpha_i}}) + \cdots + 
A_{N}(P_i)Y_{N}(Q_{i,{\alpha_i}}) = 0. 
\end{eqnarray*}   
We can rewrite this system of equations as 
\[\bmatrix
Y_1(Q_{i,1}) & \cdots & Y_{N}(Q_{i,1}) \\
\vdots & & \vdots \\
Y_1(Q_{i,{\alpha_i}}) & \cdots & Y_{N}(Q_{i,{\alpha_i}})
\endbmatrix
\bmatrix
A_1(P_i) \\
\vdots \\
A_{N}(P_i)
\endbmatrix = 
\bmatrix
0 \\
\vdots \\
0
\endbmatrix.
\]
The maximum 
number of linear restrictions $\pipi$ can place on the forms of 
$R_{l,\bj}$ that pass through $\X$ is simply the rank of the matrix
on the left.  By Proposition
~\ref{rank=hilbertforpnk} the rank of this matrix
is  equal to $H_{\qpi}(\bj)$.  \hfill{$\diamondsuit$}

By the claim, for each $P_i \in \pi_1(\X)$, the set
$\pi_1^{-1}(P_i)$ imposes
at most $H_{Q_{P_i}}(\bj)$ 
linear restrictions on the elements of $R_{l,\bj}$ that pass through
$\X$.  Hence the set $\X$ imposes at most ${\displaystyle 
\sum_{P_i \in \pix} H_{Q_{P_i}}(\bj)}$ linear restrictions.  
We thus have
\[\dim_{\bf k} (I_{\X})_{l,\bj} \geq \dim_{\bf k} R_{l,\bj} - (*),\]
or equivalently, $\dim_{\bf k} (R/I_{\X})_{l,\bj} \leq (*)$ for all integers $l$.

We will now show that if $l = t_1-1$, then the bound $(*)$ is attained.  
The set $\pix= \{P_1,\ldots,P_{t_1}\}$ is a subset of $\pr^{n_1}$.  
By Corollary ~\ref{seperators}, 
there exist $t_1$ forms $F_{P_1},\ldots,F_{P_{t_1}}$
of degree $t_1-1$ in 
${\bf k}[x_0,\ldots,x_{n_1}]$ such that $F_{P_i}(P_i) \neq 0$
and $F_{P_i}(P_j) = 0$ if $i \neq j$.  Under the natural 
inclusion ${\bf k}[x_0,\ldots,x_{n_1}] \hookrightarrow R$
we can consider the forms 
$F_{P_1},\ldots,F_{P_{t_1}}$ as forms
of $R$ of degree $(t_1,\underline{0})$.

For our fixed $\bj$, we partition the points of $\pix$ as follows:
\[S_h :=   \left\{P_i \in \pix \left| H_{\qpi}(\bj) = h\right\}\right. 
\hspace{.5cm}\text{for $h = 1,\ldots,N$.}\]  Pick a point $P_i \in \pix$
and suppose that $P_i \in S_h$ and suppose that $\qpi = 
\{Q_{i,1},\ldots,Q_{i,{\alpha_i}}\} \subseteq \pntwo$.  By using 
Proposition  ~\ref{seperators3}, 
there exists a subset $Q \subseteq \qpi$ of
$h$ elements, say $Q = \{Q_{i,1},\ldots,Q_{i,h}\}$
after a possible reordering, such that for every $Q_{i,d}
\in Q$ there exists a form  $G_{Q_{i,d}} \in 
{\bf k}[x_{2,0},\ldots,x_{2,n_2},\ldots,x_{k,0},\ldots,x_{k,n_k}]$
of degree $\bj$ such that $G_{Q_{i,d}}(Q_{i,d}) \neq 0$ but
$G_{Q_{i,d}}(Q_{i,e}) = 0$ if $Q_{i,e} \neq Q_{i,d}$ and $Q_{i,e} \in Q$.  
Under the natural
inclusion 
${\bf k}[x_{2,0},\ldots,x_{2,n_2},\ldots,x_{k,0},\ldots,x_{k,n_k}]
 \hookrightarrow R$
we can consider each $G_{Q_{i,d}}$ as an element of 
$R$ of degree $(0,\bj)$.
From this $P_i$ and  subset
$Q \subseteq Q_{P_i}$ we construct the set of forms
\[\mathcal{B}_{P_i} := \left\{F_{P_i}G_{Q_{i,1}},\ldots,F_{P_i}G_{Q_{i,h}}
\right\}.\]
We observe that 
$F_{P_i}G_{Q_{i,d}} \not\in I_{\X}$ for $d = 1,\ldots,h$ because
it fails to vanish at $P_i \times Q_{i,d}$.  Moreover, each element
of $\mathcal{B}_{P_i}$ has degree $(t_1,\bj)$ and 
$|\mathcal{B}_{P_i}| = H_{\qpi}(\bj) = h$.

We repeat the above construction for every $P_i \in \pix$ and let
\[\mathcal{B} := \bigcup_{P_i \in {\pix}} \mathcal{B}_{P_i}.\]

\noindent
{\it Claim. }  The elements of $\mathcal{B}$ are linearly independent
modulo $I_{\X}$.

\noindent
{\it Proof of the Claim. }  
It is enough to show that for each $F_{P_i}Q_{i,l} \in \mathcal{B}$,
the point $P_i \times Q_{i,l}$ does not vanish at $F_{P_i}Q_{i,l}$
but vanishes at all the 
other elements $F_{P_{i'}}Q_{{i'},{l'}} \in \mathcal{B}$.
But this follows immediately from our construction of the elements of $\mathcal{B}$.
\hfill{$\diamondsuit$}

\noindent
Because the elements of $\mathcal{B}$ are linearly independent elements modulo
$I_{\X}$ of degree $(t_1-1,\bj)$,
it follows that $\dim_{\bf k} (R/I_{\X})_{t_1-1,\bj} \geq |\mathcal{B}|.$
But since
\[|\mathcal{B}| = \sum_{P_i \in \pi_1(\X)} |\mathcal{B}_{P_i}| =
\sum_{P_i \in \pi_1(\X)} H_{Q_{P_i}}(\bj) = (*),\]
the claim implies that 
$\dim_{\bf k} (R/I_{\X})_{t_1-1,\bj} \geq (*)$.  Combining
this inequality with the previous inequality gives  
$\dim_{\bf k} (R/I_{\X})_{t_1-1,\bj} = (*)$.

To complete the proof, let $l \in \N$ be such that $l > t_1-1.$  Then, by
Proposition ~\ref{simpleboundsonhx},
we have $(*) = H_{\X}(t_1-1,\bj) \leq H_{\X}(l,\bj) \leq  (*)$,
as wanted.
\end{proof}

For each $i=1,\ldots,k$ we let $\pi_i:\pnk \rightarrow \pr^{n_i}$ be the 
projection morphism.  Set $t_i:= |\pi_i(\X)|$.  If we
partition $\X$ with respect to any of the other $(k-1)$ coordinates,
then a result identical to Proposition ~\ref{mainresult3} holds.  
Indeed, if $\bj = (j_1,\ldots,j_k)
\in \N^k$, and if we fix all but the $i^{th}$ coordinate of $\bj$,
then for all integers $l \geq t_i -1$  
\[H_{\X}(j_1,\ldots,j_{i-1},l,j_{i+1},\ldots,j_k) =
H_{\X}(j_1,\ldots,j_{i-1},t_i-1,j_{i+1},\ldots,j_k).\]

\begin{corollary}
\label{corpnk}
Let $\X$ be a set of $s$ distinct points in $\pnk$.
Fix an $i \in \{1,\ldots,k\}$.  Let
$\pi_i(\X) =\{P_1,\ldots,P_{t_i}\}$ be the
set of $t_i \leq s$ distinct $i^{th}$ coordinates in
$\X$.  Then
\begin{enumerate}
\item[$(i)$]  for all integers $l\geq t_i -1$, $H_{\X}(le_i) =
t_i$.
\item[$(ii)$]  if $j_h \gg 0$ for all $h \neq i$ and 
$j_i \geq t_i -1$, then $H_{\X}(j_1,\ldots,j_i,\ldots,j_k) = s$.
\end{enumerate}
\end{corollary}

\begin{proof}
To prove statements $(i)$ and $(ii)$, we consider only the
case that $i=1$.  The other cases will follow similarly.

Set $Q_{P_i} := \pi_{2,\ldots,k}(\pi_1^{-1}(P_i))
\subseteq \pr^{n_2}\times \cdots \times\pr^{n_k}$ 
for every $P_i \in \pi_1(\X)$, and
let $\alpha_i = |Q_{P_i}|$.  
For all sets $Q_{P_i}$, 
$H_{Q_{P_i}}(\underline{0}) = 1$.
The conclusion of
$(i)$ will follow if we use Proposition ~\ref{mainresult3} to compute
$H_{\X}(le_1)$.

To prove $(ii)$ we observe that by induction on $k$
and Proposition ~\ref{mainresult3}, if $j_h \gg 0$ for $h \neq 1$,
then $H_{Q_{P_i}}(j_2,,\ldots,j_k) = |Q_{P_i}|
= \alpha_i$ for every $P_i \in \pi_1(\X)$.  Since 
${\displaystyle \sum_{i=1}^{t_1} \alpha_i =s}$, the result is now a 
consequence of Proposition ~\ref{mainresult3}. 
\end{proof}

\begin{remark}
By  Corollary ~\ref{corpnk} $(ii)$ we have 
$H_{\X}(\bj) = s$ for all $\bj \geq (t_1-1,\ldots,t_k-1)$. 
\end{remark}

If 
$\underline{j} = (j_1,\ldots,j_k) \in \N^k$, then we denote
the vector  $(j_1,\ldots,\hat{j_i},\ldots,j_k) \in \N^{k-1}$
by $\underline{j}_i$.  Using
this notation, we have the following consequence
of Proposition ~\ref{mainresult3}.

\begin{corollary}
\label{bordercor}
Let $\X$ be a set of $s$ distinct points in $\pnk$ and
let  $t_i = |\pi_i(\X)|$ for $1 \leq i \leq k$.  
Define
$l_i := (t_1-1,\ldots,\widehat{t_i-1},\ldots,t_k-1)$ 
for $i = 1,\ldots,k$.
Then
\[H_{\X}(\bj) = \left\{
\begin{tabular}{ll}
$s$ &  $(j_1,\ldots,j_k) \geq (t_1-1,t_2-1,\ldots,t_k-1)$ \\
$H_{\X}(t_1-1,j_{2},\ldots,j_k)$ & if $j_1 \geq t_1-1$ and 
$\underline{j}_1 \not\geq l_1$ \\
~~$\vdots$ & ~~$\vdots$\\
$H_{\X}(j_1,\ldots,j_{i-1},t_i-1,j_{i+1},\ldots,j_k)$ &
  if $j_i \geq t_i-1$ and 
$\underline{j}_i \not\geq l_i$  \\
~~$\vdots$ & ~~$\vdots$\\
$H_{\X}(j_1,\ldots,j_{k-1},t_k-1)$ & if $j_k \geq t_k-1$ and 
$\underline{j}_k \not\geq l_k$ \\
\end{tabular}
\right..
\]
\end{corollary}

\begin{remark}
\label{remark_preborder}
Suppose $\X$ is a set of distinct points in $\pnk$.  Let
$\bj = (j_1,\ldots,j_k) \in \N^k$ and suppose
that $j_1 \geq t_1-1 = |\pi_1(\X)|-1$, and $j_2 \geq
t_2 -1 = |\pi_2(\X)|-1.$  Then Corollary ~\ref{bordercor} 
implies that
\[H_{\X}(j_1,j_2,\ldots,j_k)  =  H_{\X}(t_1-1,j_2,\ldots,j_k) 
 =  H_{\X}(t_1-1,t_2-1,j_3,\ldots,j_k).\]
More generally, to compute $H_{\X}(\bj)$, the above corollary
implies that if $j_i \geq t_i -1$, then we can
replace $j_i$ with $t_i -1$ and compute the Hilbert function
at the resulting tuple.  Therefore, to completely determine
$H_{\X}$ for all $\bj \in \N^k$, we need to compute
the Hilbert function only for  $\bj \leq (t_1-1,\ldots,t_k-1)$.  Since
there are only ${\displaystyle \left(\prod_{i=1}^{k} t_i\right)}$ 
$k$-tuples in $\N^k$ that have this
property, we therefore need to compute only a finite number
of values.  Furthermore, since $t_i = |\pi_i(\X)|$, the $k$-tuples
of $\N^k$ for which we need to compute the Hilbert function is determined 
from crude numerical information about $\X$,
namely the sizes of the sets $\pi_i(\X)$.  Hence, Corollary
~\ref{bordercor} is the desired generalization of Proposition
~\ref{eventual} to points in $\pnk$.
\end{remark}

For the moment, we will specialize to the case that $\X$
is a set of distinct points in $\pnpm$.  
In this context, the above corollary becomes
\[H_{\X}(i,j) = \left\{ 
\begin{tabular}{cl}
$s$ & if $(i,j) \geq (t-1,r-1)$ \\
$H_{\X}(t-1,j)$ & if $i \geq t-1$ and $j < r-1$ \\
$H_{\X}(i,r-1)$ & if $j \geq r-1$ and $i < t-1$ \\
\end{tabular}
\right.
\]
where $t = |\pix|$ and $r = |\pixt|$.
Thus, it follows that if we know $H_{\X}(t-1,j)$ for
$j = 0,\ldots,r-1$ and $H_{\X}(i,r-1)$ for 
$i = 0, \ldots,t-1$, then we know the Hilbert function
for all but a finite number of $(i,j) \in \N^2$.  This 
observation motivates the next definition.

\begin{definition}
\label{borderpnpm}
Suppose $\X \subseteq \pnpm$ is a set of $s$ distinct points and let 
$t = |\pix|$ and $r = |\pixt|$.  
Suppose that $H_{\X}$ is the Hilbert function
of $\X$.
We call the tuples
\[B_C  :=  (H_{\X}(t-1,0),H_{\X}(t-1,1),\ldots, H_{\X}(t-1,r-1))\]
and
\[B_R  :=  (H_{\X}(0,r-1),H_{\X}(1,r-1),\ldots,H_{\X}(t-1,r-1))\]
the {\it eventual column vector} and {\it eventual row vector}
respectively.  Let $B_{\X} := (B_C,B_R)$.
We call $B_{\X}$ the {\it border} of the Hilbert function
of $\X \subseteq \pnpm$. 
\end{definition}

The term border is inspired by the ``picture'' of $H_{\X}$ if
we visualize $H_{\X}$ as an infinite matrix $(m_{i,j})$ where $m_{i,j} =
H_{\X}(i,j)$.  Indeed, if
$\X \subseteq \pnpm$ with $|\pix| = t$ and $|\pixt| = r$, then
\[H_{\X} = \bmatrix
&  &   & {\bf m_{0,r-1}} & m_{0,r-1} & \cdots \\
&*  &   & {\bf m_{1,r-1}} & m_{1,r-1} & \cdots \\
&  &   & \vdots & \vdots & \\
{\bf m_{t-1,0}} & {\bf m_{t-1,1}} & \cdots  & {\bf m_{t-1,r-1} = s } 
& s & \cdots \\
m_{t-1,0} & m_{t-1,1} & \cdots  & s  & s & \cdots \\
\vdots & \vdots & & \vdots&\vdots&\ddots
\endbmatrix
.\]
The bold numbers form the border $B_{\X}$.  The entries $m_{i,j}$
with $(i,j) \leq (t-1,r-1)$ are either ``inside'' the border or
entries of the border, 
and need to be determined.  Entries with $(i,j) \geq (t,0)$
or $(i,j) \geq (0,r)$ are ``outside'' the border.  These
values depend only upon the values in the border $B_{\X}$.

The term eventual column vector is given to $B_C = (m_{t-1,0},\ldots,
m_{t-1,r-1})$ because the $i^{th}$ entry of $B_C$ is the value at
which  the $(i-1)^{th}$ column stabilizes (because our
indexing starts at zero).  We christen $B_R$
the eventual row  vector to capture a similar result about the rows.
From Corollary ~\ref{corpnk} we always have
\[B_C  =  (t,m_{t-1,1},\ldots, m_{t-1,r-2},s)
\hspace{.5cm}\text{and}\hspace{.5cm}
B_R  =  (r,m_{1,r-1},\ldots,m_{t-2,r-1},s).\]

Suppose now that $\X$ is a set of distinct points in $\pnk$.
By Remark  ~\ref{remark_preborder} it follows that if we know 
the values of $H_{\X}(t_1-1,j_2,\ldots,j_k)$ for
all $(j_2,\ldots,j_k) \leq (t_2 -1,\ldots,t_k-1)$, 
$H_{\X}(j_1,t_2-1,\ldots,j_k)$ for
all $(j_1,j_3,\ldots,j_k) \leq (t_1-1,t_3 -1,\ldots,t_k-1),\ldots,$ and 
$H_{\X}(j_1,\ldots,j_{k-1},t_k-1)$ for
all $(j_1,\ldots,j_{k-1}) \leq (t_1-1,\ldots,t_{k-1}-1)$,
then we know the Hilbert function of $\X$ for all but a finite
number of $\bj \in \N^k$.  From this observation we can extend the definition
of a border to the Hilbert functions of sets of points in $\pnk$.

\begin{definition}
\label{borderforpnk}
Let $\X$ be a set of $s$ distinct points
in $\pnk$, and let
$t_i = |\pi_i(\X)|$ for
$i = 1,\ldots,k$.  Suppose that  $H_{\X}$ is the Hilbert function of $\X$.
For each $1 \leq i \leq k$, let $B_i := \left(b_{j_1,\ldots,\hat{j_i},\ldots,
j_k}\right)$ be the $(k-1)$-dimensional array of size
$t_1 \times \cdots \times \hat{t_i} \times \cdots t_k$ where
\[b_{j_1,\ldots,\hat{j_i},\ldots,j_k} := H_{\X}(j_1,\ldots,j_{i-1},t_i-1,
j_{i+1},\ldots,j_k) ~\text{with}~ 0 \leq j_h \leq t_h -1.\]
We call $B_i$ the {\it $i^{th}$ border array} of the Hilbert function
of $\X$.  We define $B_{\X} := (B_1,\ldots,B_k)$ to be the
{\it border} of the Hilbert function of $\X$.
\end{definition}
 
\begin{remark}
If $k=2$, then $B_1$ and $B_2$ are $1$-dimensional arrays, i.e., vectors.
It is a simple exercise to verify that $B_1$ is equal
to the eventual column vector $B_C$, and $B_2 = B_R$, 
the eventual row vector, as defined in Definition ~\ref{borderpnpm}.
\end{remark}

A natural question about the entries in the border arises:

\begin{question}
\label{borderquestion}
Suppose $B = (B_1,\ldots,B_k)$ is a tuple where each
$B_i$ is a $(k-1)$-dimensional array.  Under what conditions
is $B$ the border of the Hilbert function of a set of points
in $\pnk$?
\end{question}

We would like to classify those tuples that arise as the border of a 
set of points in $\pnk$.  An answer to the above question would impose a 
severe restriction on what could be the Hilbert function of a set of points.  
This question, although weaker, is still difficult.   In the next section, 
we answer Question ~\ref{borderquestion} for the case of points in 
$\popo$.  In general, however, this problem remains open.


\section{The Border of Points in $\popo$}

Sets of points in $\popo$ and their Hilbert function were first
investigated by Giuffrida, {\it et al.} \cite{GuMaRa1}.
In this section we examine the border of the Hilbert
functions of sets of points in
$\popo$.  We first show that the border of the Hilbert function
of a set of distinct points in $\popo$, and thus, all but
a finite number of values of the Hilbert function, can be computed
directly from numerical information describing the set of points.  
We also answer Question ~\ref{borderquestion} for
sets of distinct points in $\popo$.  In this context, Question 
~\ref{borderquestion} specializes to:

\begin{question}
\label{borderquestionp1p1}
Suppose $B = (B_C,B_R)$ is a tuple where $B_C$ and $B_R$
are two vectors.  Under what conditions
is $B$ the border of the Hilbert function of a set of points
in $\popo$?
\end{question}

To demonstrate these results, we are required to recall some
definitions and results about partitions and $(0,1)$-matrices.  Once
we have recalled the relevant information, we will answer
Question ~\ref{borderquestionp1p1}.

\subsection{Partitions and $(0,1)$-matrices}
This purpose of this section is to acquaint the reader
with some results from combinatorics.   Our main reference is Ryser \cite{R}.

\begin{definition}
\label{partitiondefinition}
A tuple $\lambda = (\lambda_1,\ldots,\lambda_r)$ of positive integers
is a {\it partition} 
of an integer $s$ if $\sum \lambda_i = s$ and $\lambda_i \geq
\lambda_{i+1}$ for every $i$.  We write $\lambda = 
(\lambda_1,\ldots,\lambda_r) \vdash s$.  The {\it conjugate} of $\lambda$
is the tuple $\lambda^* = (\lambda^*_1,\ldots,\lambda^*_{\lambda_1})$
where $\lambda_i^* = \#\{\lambda_j \in \lambda ~|~ \lambda_j \geq i\}$.
Furthermore, $\lambda^* \vdash s$. 
\end{definition}

\begin{definition}
\label{ferrers}
To any partition $\lambda = (\lambda_1,\ldots,\lambda_r) \vdash s$
we can associate the following diagram:  on an $r \times \lambda_1$
grid, place $\lambda_1$ points on the first line,
$\lambda_2$ points on the second, and so on.  The resulting diagram is 
called the {\it Ferrers diagram} of $\lambda$.
\end{definition}

\begin{example}
Suppose $\lambda = (4,4,3,1) \vdash 12$.  Then the Ferrers diagram
is
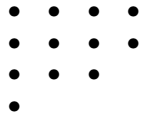
\[
\begin{tabular}{cccc}
$\bullet$ & $\bullet$ & $\bullet$ & $\bullet$ \\
$\bullet$ & $\bullet$ & $\bullet$ & $\bullet$ \\
$\bullet$ & $\bullet$ & $\bullet$ & \\
$\bullet$ & & &
\end{tabular}\]
The conjugate of $\lambda$ can be read off the Ferrers diagram by
counting the number of dots in each column as opposed to each row.  In
this example $\lambda^* = (4,3,3,2)$.
\end{example}

\begin{definition}
\label{majorizes}
Let $\lambda = (\lambda_1,\ldots,\lambda_t)$ and $\delta = 
(\delta_1,\ldots,\delta_r)$ be two partitions of $s$.  If one
partition is longer, we add zeroes to the shorter one until they
have the same length.  We say $\lambda$ {\it majorizes} $\delta$,
written $\lambda \unrhd \delta$, if 
\[ \lambda_1 + \cdots + \lambda_i \geq \delta_1 + \cdots + \delta_i 
\mbox{ for $i = 1, \ldots,\max\{t,r\}$}.\]
Majorization induces a partial ordering on the set
of all partitions of $s$.
\end{definition}

\begin{definition}
A matrix $A$ of size $m \times n$ is a {\it $(0,1)$-matrix} if all
of its entries are either zero or one.  The sum of the entries in
column $j$
will be denoted by $\alpha_j$, and the sum of the entries
of row $i$ will be denoted by
$\beta_i$.  We call the vector $\alpha_A = (\alpha_1,\ldots,\alpha_n)$
the {\it column sum vector} and the vector
$\beta_A = (\beta_1,\ldots,\beta_m)$ the {\it row sum vector}.
\end{definition}

Given a $(0,1)$-matrix, we can rearrange the rows and columns so
that $\alpha_A$ (respectively, $\beta_A$) has the property
$\alpha_i \geq \alpha_{i+1}$ (respectively $\beta_i \geq \beta_{i+1}$)
for every $i$.  Observe that $\alpha_A$ and $\beta_A$ are partitions
of the number of $1$'s in $A$.  Unless otherwise specified, we assume
that any $(0,1)$-matrix has been rearranged into this form.

If $\alpha$ and $\beta$ are any two partitions of $s$, then
we define 
\[\mathcal{M}(\alpha,\beta) := \{ (0,1)\operatorname{-matrices} ~A ~|~
\alpha_A = \alpha,  \beta_A = \beta\}.\]  
It is not evident that such 
a set is nonempty.  The following result is a classical
result, due to Gale and Ryser,  that gives us
a criterion to determine if $\mathcal{M}(\alpha,\beta) = \emptyset$.

\begin{theorem} {\em (Gale-Ryser Theorem)}
\label{ryser}
Let $\alpha$ and $\beta$ be two partitions of $s$.  The class
$\mathcal{M}(\alpha,\beta)$ is nonempty if and only if 
$\alpha^* \unrhd \beta$.
\end{theorem}
\begin{proof}
See Theorem 1.1 in Chapter 6 of Ryser's book ~\cite{R}.
\end{proof}

The proof given by Ryser to demonstrate that $\alpha^* \unrhd \beta$
implies $\mathcal{M}(\alpha,\beta)$ is nonempty is a constructive
proof.  We illustrate this construction with an example.

\begin{example}  
\label{making01matrix}
Let $\alpha =(3,3,2,1)$ and $\beta = (3,3,1,1,1)$.
A routine check will show that $\alpha^* = (4,3,2) \unrhd
(3,3,1,1,1) = \beta$.  We construct
a $(0,1)$-matrix with column sum vector $\alpha$ and row sum vector $\beta$.
Let $M$ be an empty $|\beta| \times |\alpha| = 5 \times 4$ matrix.  On
top of the $j^{th}$ column place the integer $\alpha_j$.  Beside
the $i^{th}$ row, place $\beta_i$ $1$'s.  For our example we have
\[\begin{array}{ll}
&\begin{array}{llllll}
& 3 & 3 & 2 & 1&
 \end{array}  \\

\begin{array}{l}
1~ 1 ~1 \\
1~ 1~ 1 \\
1 \\
1 \\
1 
\end{array} 
 &
\left(\begin{array}{llllll}
& & & & & \\
& & & & &\\
& & & & &\\
& & & & & \\
\end{array}\right)
\end{array}.\]
Starting with the rightmost column, we see that this column needs
one $1$.  Move a $1$ from the row with the largest number of $1$'s
to this column and fill the remainder of the column with zeroes.  
If two rows have the same number of ones, we take the first such row.  So, 
after one step,
\[\begin{array}{ll}
&\begin{array}{llllll}
& 3 & 3 & 2 & ~1&
 \end{array}  \\
\begin{array}{l}
1~ 1 \not1 \\
1~ 1 ~1 \\
1 \\
1 \\
1 
\end{array} 
 &
\left(\begin{array}{llllll}
&& & & &1   \\
&& & & &0 \\
&& & & &0 \\
&& & & &0  \\
&& & & &0  \\
\end{array}\right)
\end{array}.\]
We now repeat the above procedure on the next to last column.  We place
two $1$'s in the third column, taking our $1$'s from
the rows that contain the largest number of ones.  Thus, our
example becomes
\[\begin{array}{ll}
&\begin{array}{llllll}
& 3 & ~3 & ~2 & ~1&
 \end{array}  \\

\begin{array}{l}
1 \not1 \not1 \\
1~ 1 \not1 \\
1 \\
1 \\
1 
\end{array} 
 &
\left(\begin{array}{llllll}
&& & &1 &1   \\
&& & &1 &0 \\
&& & &0 &0 \\
&& & &0 &0  \\
&& & &0 &0  \\
\end{array}\right)
\end{array}.\]
We continue the above method for the remaining columns to get
\[\begin{array}{ll}
&\begin{array}{llllll}
& 3 & 3 & 2 & 1&
 \end{array}  \\

\begin{array}{l}
\not1 \not1 \not1 \\
\not1 \not1 \not1 \\
\not1 \\
\not1 \\
\not1 
\end{array} 
 &
\left(\begin{array}{llll}
0 &1 &1 &1   \\
1 &1 &1 &0 \\
0 &1 &0 &0 \\
1 &0 &0 &0  \\
1 &0 &0 &0  \\
\end{array}\right)
\end{array}.\]
It follows immediately that our matrix is an element of $\mathcal{M}(\alpha,
\beta)$.  The proof of the Gale-Ryser Theorem
shows that if $\alpha^* \unrhd \beta$,
then this algorithm always works.
\end{example}

\subsection{Classifying the borders of Hilbert Functions
of Points in $\popo$}

In Section 4 we defined the border of a Hilbert
function for points $\X \subseteq \pnk$.  Question ~\ref{borderquestion}
asks what tuples can be the border of a  Hilbert function
of a set of points.  For points $\X \subseteq \pnpm$ this question
reduces to describing all possible eventual column vectors $B_C$
and eventual row vectors $B_R$.  We wish to answer this question
for points in $\popo$.

So, suppose that $\X \subseteq \popo$ is a collection of $s$
distinct points.  We associate to $\X$ two tuples, $\ax$ and $\bx$,
as follows. For each $P_i \in \pix =\{P_1,\ldots,P_t\}$
we set $\alpha_i := |\pipi|$.  After relabeling the $\alpha_i$'s so
that $\alpha_i \geq \alpha_{i+1}$ for $i=1,\ldots,t-1$,  we set
$\ax := (\alpha_1,\ldots,\alpha_t)$.  Analogously, for every
$Q_i \in \pixt = \{Q_1,\ldots, Q_r\}$ we set $\beta_i := |\piqi|$.
After relabeling the $\beta_i$'s so that $\beta_i \geq \beta_{i+1}$
for $i = 1,\ldots,r-1$, we let $\bx$ be the $r$-tuple $\bx :=
(\beta_1,\ldots,\beta_r)$. We note that 
$\ax$ and $\bx$ are both partitions (see Definition ~\ref{partitiondefinition})
of the integer $s = |\X|$.  Thus, we can write 
$\ax \vdash s$ and $\bx \vdash s$.  
If we denote the length of $\ax$ (resp. $\bx$) by
$|\ax|$ (resp. $|\bx|$), then
we also observe that
$|\pix| = |\ax|$ and $|\pixt| = |\bx|$.

\begin{remark} 
\label{altsum}
The following observation will be useful
about a set of points $\X \subseteq \pnpm$.
Fix an integer $j \geq 0$, let $\pix = \{P_1,\ldots,P_t\}$, and let 
${\displaystyle (*) = \sum_{P_i \in \pi_1(\X)} H_{\qpi}(j)}$.
It is sometimes useful to note that $(*)$ is equal to
\begin{eqnarray*}
(*) & = & \#\left\{P_i \in \pix ~\left|~ H_{\qpi}(j) =1\right\}\right. + \cdots\\
& & +
\binom{m+j}{j}\#\left\{P_i \in\pix ~\left|~ H_{\qpi}(j) =\binom{m+j}{j}
\right\}\right.,  
\end{eqnarray*}
and that $(*)$ is also equal to
\begin{eqnarray*}
\#\left\{P_i \in\pix ~\left| H_{\qpi}(j) \geq 1\right\}\right. & + & 
\#\left\{P_i \in\pix ~\left| H_{\qpi}(j) \geq 2\right\}\right. 
+\cdots \\
&+ &\#\left\{P_i \in\pix ~\left| H_{\qpi}(j) \geq \binom{m+j}{j}\right\}\right..
\end{eqnarray*}
\end{remark}

As an application of Proposition ~\ref{mainresult3}
we demonstrate that for points $\X \subseteq \popo$ the eventual
column vector $B_C$ and the eventual row vector $B_R$ can be computed
directly from the tuples $\ax$ and $\bx$.  We first recall
the Hilbert function for sets of points in $\pr^1$.

\begin{lemma}[\cite{GeMa} ~Proposition 1.3] 
\label{pointsinp1}
Let $\X = \{P_1,\ldots,P_s\}\subseteq \pr^1$.  Then
\[ H_{\X}(i) = \left\{\begin{array}{ll}
i+1 & 0 \leq i \leq s-1 \\
s & i \geq s
\end{array}
\right..\]
\end{lemma}

\begin{proposition}
\label{computingborderinp1p1}
Let $\X \subseteq \popo$ be a set of $s$ distinct points and suppose
that  $\ax =
(\alpha_1,\ldots,\alpha_t)$ and $\bx = (\beta_1,\ldots,\beta_r)$.  
Let $B_C = (b_0,b_1,\ldots,b_{r-1})$ where $b_j = H_{\X}(t-1,j)$,
be the eventual column vector of the Hilbert function $H_{\X}$.  Then
\[b_j = \#\{\alpha_i \in \ax ~|~ \alpha_i \geq 1\} +  
\#\{\alpha_i \in \ax ~|~ \alpha_i \geq 2\} + \cdots +  
\#\{\alpha_i \in \ax ~|~ \alpha_i \geq j+1\}.\]
Analogously, if $B_R = (b'_0,b'_1,\ldots,b'_{t-1})$, with 
$b'_j = H_{\X}(j,r-1)$, is the eventual row vector of $H_{\X}$, then
\[b'_j = \#\{\beta_i \in \bx ~|~ \beta_i \geq 1\} +  
\#\{\beta_i \in \bx ~|~ \beta_i \geq 2\} + \cdots +  \#\{\beta_i \in \bx 
~|~ \beta_i \geq j+1\}.\]
\end{proposition}
\begin{proof}
After relabeling the elements of $\pi_1(\X)$, we can
assume that $|\pi_1^{-1}(P_i)| = \alpha_i$.
By Proposition ~\ref{mainresult3} and Remark ~\ref{altsum} we have
\begin{eqnarray*}
b_j = H_{\X}(t-1,j) & = &
\#\left\{P_i \in \pix ~ \left| H_{\qpi}(j) \geq 1\right\}\right. + 
\#\left\{P_i \in \pix ~\left| H_{\qpi}(j) \geq 2\right\}\right. 
\\
& & 
+ \cdots + \#\left\{P_i \in \pix ~\left| H_{\qpi}(j) \geq j+1\right\}\right..
\end{eqnarray*}
Now $Q_{P_i} = \pi_2(\pipi)$ is a subset of $\alpha_i$ points in
$\pr^1$.  If $1 \leq k \leq j+1$, then $H_{Q_{P_i}}(j) \geq k$ if
and only if $|\pi_1^{-1}(P_i)| \geq k$.  This is a consequence of Lemma
~\ref{pointsinp1}.  This in turn implies that the
sets 
\[\left\{P_i \in \pix~\left|~ H_{\qpi}(j) \geq k\right\}\right. 
\hspace{.5cm}\text{and}\hspace{.5cm} 
\left\{P_i \in \pix ~\left|~ |\pi_1^{-1}(P_i)| \geq k\right\}\right.\] 
are the same, and thus, the numbers
$\#\{P_i \in \pix~|~ H_{\qpi}(j) \geq k\}$ and 
$\#\{\alpha_i \in \ax ~|~ \alpha_i \geq k\}$ are equal.  The desired
identity now follows from this result.  The statement about
the eventual row vector $B_R$ is proved similarly.
\end{proof}

We can rewrite the above result more succinctly by invoking the language
of combinatorics introduced earlier.  Recall that the conjugate of a partition
$\lambda = (\lambda_1,\ldots,\lambda_k)$ is the tuple
$\lambda^* = (\lambda^*_1,\ldots,\lambda^*_{\lambda_1})$ where
$\lambda^*_j := \#\{\lambda_i \in \lambda ~|~ \lambda_i \geq j\}$.

\begin{definition}
If $p = (p_1,p_2,\ldots,p_k)$, then 
$\Delta p := (p_1,p_2 - p_1, \ldots, p_k - p_{k-1}).$
\end{definition}

\begin{corollary}
\label{borderinp1p1}
Let $\X \subseteq \popo$ be $s$ distinct points with $\ax$ and $\bx$.
Then
\begin{enumerate}
\item[$(i)$]  $\Delta B_C = \ax^*$.
\item[$(ii)$]  $\Delta B_R = \bx^*$.
\end{enumerate}
\end{corollary}
\begin{proof}
Using Proposition 
~\ref{computingborderinp1p1} to calculate $\Delta B_C$ we get
\[\Delta B_C = ( \#\{\alpha_i \in \ax ~|~ \alpha_i \geq 1\},
\#\{\alpha_i \in \ax ~|~ \alpha_i \geq 2\},\ldots, 
\#\{\alpha_i \in \ax ~|~ \alpha_i \geq r\}).\]
The conclusion follows by noting that
 $ \#\{\alpha_i \in \ax ~|~ \alpha_i \geq j\}$ is by definition
the $j^{th}$ coordinate of $\ax^*$.
The proof of $(ii)$ is the same as $(i)$. 
\end{proof}

\begin{remark}
\label{corborderremark}
For each positive integer $j$ we have the
following identity:
\[
\#\{\alpha_i \in \ax ~|~ \alpha_i \geq j\} - 
\#\{\alpha_i \in \ax ~|~ \alpha_i \geq j+1\}
=
\#\{\alpha_i \in \ax ~|~ \alpha_i =j\}
.\]
Since Corollary ~\ref{borderinp1p1} shows that
\[\#\{\alpha_i \in \ax ~|~ \alpha_i \geq j\} =
 H_{\X}(t-1,j-1)-H_{\X}(t-1,j-2)\]
it follows from the above identity that 
\begin{eqnarray*}
\#\{\alpha_i \in \ax ~|~ \alpha_i = j\}
& = & 
\left[H_{\X}(t-1,j-1)-H_{\X}(t-1,j-2)\right]- \\
& &
\left[H_{\X}(t-1,j) - H_{\X}(t-1,j-1)\right].
\end{eqnarray*}
Thus, for each integer $1 \leq j \leq r$ there is precisely
$[H_{\X}(t-1,j-1)-H_{\X}(t-1,j-2)] - [H_{\X}(t-1,j) - H_{\X}(t-1,j-1)]$
lines of degree $(1,0)$
that pass through $\X$ that contain exactly $j$ points of $\X$.
This is the statement of Theorem 2.12 of Giuffrida, {\it et al.}
~\cite{GuMaRa1}.  Of course, a similar result holds for
the lines of degree $(0,1)$ that pass through $\X$.
\end{remark}

\begin{example}
\label{nonacmexample}
We illustrate how to use Corollary ~\ref{borderinp1p1}
 to compute the Hilbert function
for a set of points $\X \subseteq \popo$ for all
but a finite number $(i,j) \in \N^2$.  Suppose that 

\begin{center}
\begin{picture}(150,100)(25,0)
\put(10,40){$\X = $}
\put(60,10){\line(0,1){70}} 
\put(80,10){\line(0,1){70}}
\put(100,10){\line(0,1){70}}
\put(120,10){\line(0,1){70}}
\put(140,10){\line(0,1){70}}
\put(160,10){\line(0,1){70}}

\put(54,-5){$P_1$}
\put(74,-5){$P_2$}
\put(94,-5){$P_3$}
\put(114,-5){$P_4$}
\put(134,-5){$P_5$}
\put(154,-5){$P_6$}

\put(55,15){\line(1,0){115}}
\put(55,35){\line(1,0){115}}
\put(55,55){\line(1,0){115}}
\put(55,75){\line(1,0){115}}

\put(180,11){$Q_1$}
\put(180,31){$Q_2$}
\put(180,51){$Q_3$}
\put(180,71){$Q_4$}

\put(60,15){\circle*{5}}
\put(80,15){\circle*{5}}
\put(80,35){\circle*{5}}
\put(80,55){\circle*{5}}
\put(80,75){\circle*{5}}
\put(100,35){\circle*{5}}
\put(120,15){\circle*{5}}
\put(120,55){\circle*{5}}
\put(140,55){\circle*{5}}
\put(140,75){\circle*{5}}
\put(160,35){\circle*{5}}
\put(160,55){\circle*{5}}
\put(160,75){\circle*{5}}
\end{picture}
\end{center}

\noindent
For this example $\alpha_{\X} = (4,3,2,2,1,1)$ because
$|\pi_1^{-1}(P_1)| = 1$, 
$|\pi_1^{-1}(P_2)| = 4$,
$|\pi_1^{-1}(P_3)| = 1$,
$|\pi_1^{-1}(P_4)| = 2$,
$|\pi_1^{-1}(P_5)| = 2$, and 
$|\pi_1^{-1}(P_6)| = 3$.  The conjugate of $\ax$ is $\ax^* = (6,4,2,1)$,
and hence, by Corollary ~\ref{borderinp1p1} we know that
$B_C = (6,10,12,13).$  Similarly, $\bx = (4,3,3,3)$, and thus
$\bx^* = (4,4,4,1)$.  Using Corollary ~\ref{borderinp1p1} we have
$B_R = (4,8,12,13,13,13)$.  (Note that we need to add some $13$'s to the end
of $B_R$ to ensure that $B_R$ has the correct length of 
$|B_R| = |\pix| = 6$.)  Visualizing the Hilbert function $H_{\X}$
as a matrix and using the tuples $B_R$ and $B_C$, we have
\[H_{\X} = \bmatrix
&  &   & {\bf 4} & 4 & \cdots \\
& * &   & {\bf 8} & 8 & \cdots \\
&  &   & {\bf 12} & 12 & \cdots\\
&  &   & {\bf 13} & 13 & \cdots\\
&  &   & {\bf 13} & 13 & \cdots\\
{\bf 6}& {\bf 10} & {\bf 12}   & {\bf 13} & 13 & \cdots\\
 6&  10 & 12   &  13 & 13 & \cdots\\
\vdots& \vdots  & \vdots  & \vdots & \vdots &\ddots\\
\endbmatrix
.\]
All that remains to be calculated are the entries in the upper left-hand
corner of $H_{\X}$ denoted by $*$.
\end{example}

As is evident from Corollary ~\ref{borderinp1p1} and 
Remark ~\ref{corborderremark}, 
the border of the Hilbert function for points $\X \subseteq \popo$ 
is linked to combinatorial information describing some of the
geometry of $\X$, e.g., the number of points whose first coordinate
is $P_1$, the number of points whose first coordinate is $P_2$,
etc.  By utilizing the Gale-Ryser Theorem (Proposition ~\ref{ryser})
we show that the geometry of $\X$ forces a condition
on $\ax$ and $\bx$.  As a corollary,  
we can answer Question ~\ref{borderquestion}
for points in $\popo$.
\begin{theorem}
\label{borderresult}
Let $\alpha,\beta \vdash s$.  Then there exists a set of points
$\X \subseteq \popo$ such that $\ax = \alpha$ and $\bx = \beta$
if and only if $\alpha^* \unrhd \beta$.
\end{theorem}
\begin{proof}
Suppose that there exists a set of points $\X$ such that 
$\ax = \alpha$ and $\bx = \beta$.  Suppose that $\pix = \{P_1,\ldots,P_t\}$
with $t = |\alpha|$. For $i = 1,\ldots,t$, 
let $L_{P_i}$ be the line in $\popo$
defined by the $(1,0)$-form such that $\pipi \subseteq L_{P_i}$.  Similarly,
if $\pixt = \{Q_1,\ldots,Q_r\}$, where $r = |\beta|$, let 
$L_{Q_i}$ be the line defined by the $(0,1)$-form such
that $\piqi \subseteq L_{Q_i}$.  For each pair $(i,j)$ where $1\leq i\leq t$
and $1\leq j \leq r$, the lines $L_{P_i}$ and $L_{Q_j}$
intersect at a unique point $P_i \times Q_j$.  We note
that $\X \subseteq \{ P_i \times Q_j ~|~  1 \leq i \leq t, 1\leq j \leq r\}$.
We define an $r\times t$ $(0,1)$-matrix $A=(a_{i,j})$ where
\[ a_{i,j} = \left\{ 
\begin{tabular}{ll}
1 & if $L_{P_i} \cap L_{Q_j} = P_i \times Q_j \in \X$ \\
0 &  if $L_{P_i} \cap L_{Q_j} = P_i \times Q_j \not\in \X$ 
\end{tabular}
\right..\]
By construction this $(0,1)$-matrix has column sum vector $\alpha_A = \ax$
and row sum vector $\beta_A = \bx$.  Hence, $\mathcal{M}(\alpha,\beta) \neq
\emptyset$ because $A \in \mathcal{M}(\alpha,\beta)$. 
The conclusion  $\alpha^* \unrhd \beta$ follows from  the  
Gale-Ryser Theorem (Proposition ~\ref{ryser}).

To prove the converse, it is sufficient to construct a set
$\X \subseteq \popo$ with $\ax =\alpha$ and $\bx = \beta$.  Since
$\alpha^* \unrhd \beta$ there exists a $(0,1)$-matrix $A \in 
\mathcal{M}(\alpha,\beta)$.  Fix such a matrix $A$.  Let $L_{P_1},\ldots,
L_{P_t}$ be $t = |\alpha|$ distinct lines in $\popo$
defined by forms of degree $(1,0)$, and let $L_{Q_1},\ldots,
L_{Q_r}$ be $r = |\beta|$ distinct lines in $\popo$
defined by forms of degree $(0,1)$.  For every pair $(i,j)$, with
$1 \leq i \leq t$ and $1 \leq j \leq r$, the lines $L_{P_i}$
and $L_{Q_j}$ intersect at the distinct point $P_i \times Q_j =
L_{P_i} \cap L_{Q_j}$.   We define a set of points $\X \subseteq \popo$
using the matrix $A = (a_{i,j})$ as follows:
\[ \X := \{P_i \times Q_j ~|~ a_{i,j} = 1\}.\]
From our construction of $\X$ we have $\ax = \alpha$
and $\bx = \beta$.
\end{proof}

\begin{remark}  Suppose that $\alpha,\beta \vdash s$ and $\alpha^* \unrhd
\beta$.
Then by adopting the procedure described in Example ~\ref{making01matrix},
we can construct a set of points $\X$
in $\popo$ with $\ax = \alpha$ and $\bx = \beta$.  For example,
if $\alpha = (3,3,2,1)$ and $\beta = (3,3,1,1,1)$ are as
in Example ~\ref{making01matrix}, then we saw how to construct a $(1,0)$-matrix
from $\alpha$ and $\beta$.  We then identify this matrix with
a set of points as in the proof of Theorem ~\ref{borderresult}.  For
the example of Example ~\ref{making01matrix} we have
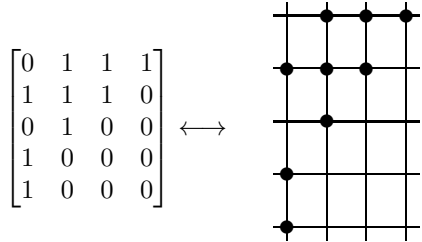
\[
\bmatrix
0 &1 &1 &1   \\
1 &1 &1 &0 \\
0 &1 &0 &0 \\
1 &0 &0 &0  \\
1 &0 &0 &0  \\
\endbmatrix
\longleftrightarrow
\begin{picture}(100,50)(0,0)

\put(20,-40){\line(0,1){90}} 
\put(35,-40){\line(0,1){90}}
\put(50,-40){\line(0,1){90}}
\put(65,-40){\line(0,1){90}}

\put(15,-35){\line(1,0){55}}
\put(15,-15){\line(1,0){55}}
\put(15,05){\line(1,0){55}}
\put(15,25){\line(1,0){55}}
\put(15,45){\line(1,0){55}}

\put(20,-35){\circle*{5}}
\put(20,-15){\circle*{5}}
\put(20,25){\circle*{5}}
\put(35,45){\circle*{5}}
\put(35,25){\circle*{5}}
\put(35,05){\circle*{5}}
\put(50,45){\circle*{5}}
\put(50,25){\circle*{5}}
\put(65,45){\circle*{5}}
\end{picture}\]
\end{remark}

\begin{corollary}
\label{borderresultcor}
Suppose $B_C = (b_0,\ldots,b_{r-1})$ and $B_R = (b'_0,\ldots,b'_{t-1})$ 
are two tuples
such that $b_0 =t$, $b'_0 = r$, and $\Delta B_C, \Delta B_R \vdash s$.  Then $B_C$ is the eventual column
vector and $B_R$ is the eventual row vector 
of a Hilbert function of a set of $s$ points in $\popo$
if and only if $\Delta B_C \unrhd (\Delta B_R)^*$.
\end{corollary}
\begin{proof}
For any partition $\lambda$, we have the identity $(\lambda^*)^* = \lambda$.
If $B_{\X} =(B_C,B_R)$ is the border of a set of points, then
$\Delta B_C = \alpha_{\X}^* \unrhd \bx = (\beta_{\X}^*)^* = (\Delta B_R)^*$. 

Conversely, suppose that $\Delta B_C \unrhd (\Delta B_R)^*$.  Let
$\alpha = (\Delta B_C)^*$ and $\beta = (\Delta B_R)^*$.  
Since $\alpha^* \unrhd \beta$, there exist a set of points
$\X \subseteq \popo$ with $\alpha_{\X} = \alpha$ and $\bx = \beta$.  But
then $\Delta B_C = \Delta B'_C$, where $B'_C$
is the eventual column vector of the Hilbert function of
$\X$.  Since $|B'_C| = |\beta| = r,$ and 
because first element of the tuple $B'_C$ is $t$, we have
$B_C = B'_C$.  We show that the eventual row border $B'_R$ of
the Hilbert function of 
$\X$ is equal to $B_R$ via the same argument.
\end{proof}

\begin{remark}
It is possible for two sets of points 
to have the same border, but not the
same Hilbert function.  For example, let $P_i :=[1:i]$ with
$i=1,2,3$ be three 
distinct points of $\pr^1$, and let $Q_i :=[1:i]$ with
$i = 1,2,3$ be another
collection of three distinct points in $\pr^1$.
Let $\X_1 = \{P_1 \times Q_1, P_2 \times Q_2, P_2 \times Q_3, 
P_3 \times Q_1\}$, and let 
$\X_2 = \{P_1 \times Q_3, P_2 \times Q_1, P_2 \times Q_2, 
P_3 \times Q_1\}$.  We can visualize these sets as 

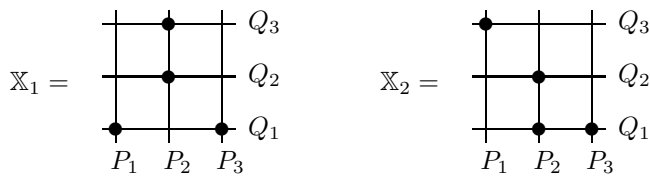
\begin{picture}(100,80)(-50,-5)
\put(0,30){$\X_1 = $}
\put(40,10){\line(0,1){50}} 
\put(60,10){\line(0,1){50}}
\put(80,10){\line(0,1){50}}
\put(38,0){$P_1$}
\put(58,0){$P_2$}
\put(78,0){$P_3$}
\put(180,10){\line(0,1){50}}
\put(200,10){\line(0,1){50}}
\put(220,10){\line(0,1){50}}
\put(178,0){$P_1$}
\put(198,0){$P_2$}
\put(218,0){$P_3$}

\put(35,15){\line(1,0){50}}
\put(35,35){\line(1,0){50}}
\put(35,55){\line(1,0){50}}
\put(90,13){$Q_1$}
\put(90,33){$Q_2$}
\put(90,53){$Q_3$}

\put(40,15){\circle*{5}}
\put(60,35){\circle*{5}}
\put(60,55){\circle*{5}}
\put(80,15){\circle*{5}}

\put(140,30){$\X_2 = $}
\put(175,15){\line(1,0){50}}
\put(175,35){\line(1,0){50}}
\put(175,55){\line(1,0){50}}
\put(230,13){$Q_1$}
\put(230,33){$Q_2$}
\put(230,53){$Q_3$}

\put(180,55){\circle*{5}}
\put(200,35){\circle*{5}}
\put(200,15){\circle*{5}}
\put(220,15){\circle*{5}}
\end{picture}

\noindent
For this example, $\alpha_{\X_1} = \alpha_{\X_2} = (2,1,1)$
and $\beta_{\X_1} = \beta_{\X_2} = (2,1,1)$, and hence, both sets of points
have the same border.  However, using 
{\tt CoCoA} to compute the Hilbert function
of $\X_1$ and $\X_2$, we find that that the
Hilbert functions are not equal.  
Specifically,
\[H_{\X_1} = \bmatrix
$1$ & $2$ & $3$ & \cdots \\
$2$ & $3$ & $4$ & \cdots\\
$3$ & $4$ & $4$ & \cdots \\
\vdots &\vdots &\vdots & \ddots
\endbmatrix
\hspace{.5cm}
H_{\X_2} = \bmatrix
$1$ & $2$ & $3$ & \cdots \\
$2$ & $4$ & $4$ & \cdots\\
$3$ & $4$ & $4$ & \cdots\\
\vdots & \vdots&\vdots & \ddots
\endbmatrix.
\]
\end{remark} 


\section*{Acknowledgments}

The results in this paper were part of the my Ph.D. thesis.  
I would like to thank D. Gregory, D. de Caen, A. Ragusa, L. Roberts,
D. Wehlau, and the audience of the Curve Seminar 
for their comments and suggestions.
Part of this work was completed at the   Universit\`a di Genova,
and I would like to thank the people there, especially
the members of the {\tt CoCoA} group,
for their hospitality.  I would especially like to thank my supervisor
Tony Geramita for introducing me to this problem and for his
encouragement and help.


\end{document}